\documentclass[10pt]{article}
\usepackage[nohead,margin=1.0in]{geometry}
\usepackage{amssymb, amsmath, amsthm,mathtools}
\usepackage{color}
\usepackage{graphicx}
\usepackage{epstopdf}
\graphicspath{{./pics/}}
\usepackage{booktabs}
\usepackage{threeparttable}
\usepackage{dsfont,verbatim}
\usepackage{url}
\usepackage{subfig}
\mathtoolsset{showonlyrefs}

\newtheorem{theorem}{Theorem}[section]
\newtheorem{lemma}[theorem]{Lemma}

\newtheorem{proposition}[theorem]{Proposition}
\newtheorem{definition}{Definition}[section]
\newtheorem{remark}{Remark}[section]
\newtheorem{example}{Example}[section]
\numberwithin{equation}{section}

\renewcommand{\d}{\mathrm{d}}
\allowdisplaybreaks
\providecommand{\abs}[1]{\left\lvert#1\right\rvert}
\providecommand{\norm}[1]{\left\lVert#1\right\rVert}

\newcommand{\R}{\mathbb{R}}
\def\phi {\varphi}

\title{Recovering Multiple Fractional Orders in Time-Fractional Diffusion in an Unknown Medium\thanks{The work of B.J. is partially supported by UK EPSRC grant EP/T000864/1, and that of Y.K. by the French National Research Agency ANR
(project MultiOnde) grant ANR-17-CE40-0029.}}

\author{Bangti Jin\thanks{Department of Computer Science, University College London, Gower Street, London WC1E 6BT, UK (b.jin@ucl.ac.uk)}\and
Yavar Kian\thanks{Aix Marseille Universit\'{e}, Universit\'{e} de Toulon, CNRS, CPT, Marseille, France (yavar.kian@univ-amu.fr)}}

\date{}

\begin{document}

\maketitle

\begin{abstract}
In this work, we investigate an inverse problem of recovering multiple orders in a time-fractional diffusion model
from the data observed at one single point on the boundary. We prove the unique recovery of the orders together with
their weights, which does not require a full knowledge of the domain or medium properties, e.g., diffusion and potential
coefficients, initial condition and source in the model. The proof is based on Laplace transform and asymptotic expansion.
Further, inspired by the analysis, we propose a numerical procedure for recovering these parameters based on a nonlinear
least-squares fitting with either fractional polynomials or rational approximations as the model function, and provide numerical experiments
to illustrate the approach for small time $t$.
\end{abstract}

\section{Introduction}
Let $\Omega\subset \mathbb{R}^d$ ($d\ge 2$) be an open bounded and connected subset with a
$C^{2\lceil\frac{d}{4}\rceil+2}$ boundary $\partial \Omega$ (the notation $\lceil r\rceil$ denotes
the smallest integer exceeding $r\in\mathbb{R}$). Consider a weak solution (in
the sense of Definition \ref{def:weak} below) $u$ of the following initial boundary value problem:
\begin{equation}\label{eqn:fde}
\left\{\begin{aligned}
\sum_{j=1}^Nr_j\partial_t^{\alpha_j}u +\mathcal{A} u &=  \sigma(t)f(x), \quad \mbox{in }\Omega\times(0,T),\\
\mathcal R u&= 0, \quad \mbox{on } \partial\Omega\times(0,T), \\
u&=u_0, \quad \mbox{in } \Omega\times \{0\}.
\end{aligned}\right.
\end{equation}
In the model, $\mathcal{A}$ is a second-order elliptic operator on the domain $\Omega$ given by
\begin{equation}\label{A}
\mathcal{A} u(x) :=-\sum_{i,j=1}^d \partial_{x_i} \left( a_{i,j}(x) \partial_{x_j} u(x) \right)+q(x)u(x),\quad  x\in\Omega,
\end{equation}
where the potential $q \in C^{2\lceil\frac{d}{4}\rceil}(\overline{\Omega})$ is strictly positive in $\overline{\Omega}$,
and the diffusion coefficient matrix $a:=(a_{i,j})_{1 \leq i,j \leq d} \in C^{1+2\lceil\frac{d}{4}\rceil}
(\overline{\Omega};\mathbb{R}^{d\times d})$ is symmetric and fulfills the following ellipticity condition
\begin{equation}\label{ell}
\exists c>0,\ \sum_{i,j=1}^d a_{i,j}(x) \xi_i \xi_j \geq c |\xi|^2,\quad \forall x \in \overline{\Omega},\ \xi=(\xi_1,\ldots,\xi_d) \in \mathbb{R}^d.
\end{equation}
For a fixed $N\in\mathbb N$ and for $j=1,\ldots,N$, we consider the constants $r_j\in(0,+\infty)$,
$0<\alpha_1<\ldots<\alpha_N<1$ and $T\in(0,+\infty)$.
In the model \eqref{eqn:fde}, the notation $\partial_t^\alpha u$ denotes the
Djrbashian-Caputo fractional derivative of order $\alpha$ in
$t$, for $\alpha \in (0,1) $, defined by (cf. \cite[p. 92]{KilbasSrivastavaTrujillo:2006}
and \cite[Section 2.3.2]{Jin:2021})
\begin{equation*}
 \partial_t^\alpha u(x,t):=\frac{1}{\Gamma(1-\alpha)}\int_0^t(t-s)^{-\alpha}\partial_su(x,s) \d s,\quad (x,t) \in \Omega
 \times(0,T),
\end{equation*}
where the notation $\Gamma(z)=\int_0^\infty s^{z-1}e^{-s}\d s$, $\Re(z)>0$, denotes Euler's Gamma function.
In addition, in the model \eqref{eqn:fde}, the notation $\mathcal R$ denotes either the Dirichlet trace $\mathcal Ru=
u_{|\partial\Omega\times(0,T)}$ or the normal derivative $\partial_{\nu_a}$ associated with the diffusion
coefficient matrix $a$ given by
\begin{equation*}
  \mathcal Ru=\partial_{\nu_a}u_{|\partial\Omega\times(0,T)}=\sum_{i,j=1}^da_{ij}\partial_{x_j}u\nu_i|_{\partial\Omega\times(0,T)},
\end{equation*}
where $\nu=(\nu_1,\ldots,\nu_d)\in\mathbb{R}^d$ denotes  the unit outward normal vector to the boundary $\partial\Omega$.
Throughout, the adjoint trace $\mathcal{R}^*$ denotes that $\mathcal R^*$ is the Dirichlet boundary trace if $\mathcal R$ correspond to
the Neumann one and $\mathcal R^*$ is the Neumann boundary trace if $\mathcal R$ correspond to the Dirichlet one.

When $N =1$, the model \eqref{eqn:fde} reduces to its single-term counterpart, i.e., with $\alpha\in (0,1)$,
\begin{equation}\label{eqn:fde0}
 \partial_t^\alpha u - \mathcal{A}u= f,\quad\mbox{in }\Omega \times (0, T ].
\end{equation}
This model has been studied extensively in the engineering, physical and mathematical literature
due to its extraordinary capability for describing anomalous diffusion phenomena \cite{MetzlerKlafter:2000}. It is the fractional analogue of
the classical diffusion equation: with $\alpha=1$, it recovers the latter, and thus inherits some of its
important analytical properties. However, it also differs considerably from the latter in the sense that, due to the
presence of the nonlocal fractional derivative term, it has limited smoothing property in space and slow
asymptotic decay in time \cite{Jin:2021}. The multi-term model \eqref{eqn:fde}
employs multiple fractional orders to improve the modeling accuracy of the single-term model
\eqref{eqn:fde0}. For example, a two-term fractional-order diffusion model was
proposed in \cite{SchumerBenson:2003} for the total concentration in solute transport, in order to describe the mobile and
immobile status of the solute. The model with two fractional derivatives
appears also naturally when describing subdiffusive motion in velocity fields \cite{MetzlerKlafter:1998}.

This work is interested in the following inverse problem for the model \eqref{eqn:fde}.
Let $u(x,t)$ be the weak solution to the model \eqref{eqn:fde} in the sense of
Definition \ref{def:weak} below. Given the observation $\mathcal{R}^*u(x_0,t)$, $t\in(0,T)$,
for some $x_0\in \partial\Omega$, can one uniquely determine the orders $\{\alpha_i\}_{i=1}^N$ and weights
$\{r_i\}_{i=1}^N$ in the model \eqref{eqn:fde}? Physically, the orders $\{\alpha_i\}$ are determined by
the inhomogeneity of the media, but it is still unclear which physical law can relate the inhomogeneity
to $\{\alpha_i\}_{i=1}^N$. Thus, in practice, one natural way is to formulate an inverse
problem of determining the parameters from the available data, e.g., $\mathcal{R}^*u(x_0, t)$, $0 < t < T$ at a monitoring point
$x_0\in \overline{\Omega}$. The short answer to the inverse problem is affirmative.
To precisely describe the results, we need a proper functional analytic framework. We define an operator
$A=\mathcal A$ acting on $L^2(\Omega)$ with its domain $D(A)$ given by
$D(A)=\{v\in H^2(\Omega):\ \mathcal A v\in L^2(\Omega),\ Rv=0\mbox{ on }\partial\Omega\}.$
Moreover, \cite[Theorem 2.5.1.1]{Grisvard:1985} implies that, for all $\ell=1,\ldots,\lceil\frac{d}{4}\rceil+1$,
\begin{equation}\label{DA}
D(A^\ell)=\{v\in H^{2\ell}(\Omega):\ \mathcal Rv=\mathcal R(\mathcal Av)=\ldots=\mathcal R(\mathcal A^{\ell-1} v)=0\}.
\end{equation}

We need the following assumption on the data in the model \eqref{eqn:fde}. The space $D(A^s)$
is defined in Section \ref{sec:prelim}.
\begin{definition}\label{def:admissible}
A tuple $(\Omega,a,q,f,u_0)$ is said to be admissible if the following conditions are fulfilled.
\begin{itemize}
\item[{\rm(i)}]  $\Omega\subset
\mathbb{R}^d$ is a $C^{2\lceil\frac{d}{4}\rceil+2}$ bounded open set, $a:=(a_{i,j})_{1 \leq i,j \leq d} \in C^{1+2\lceil\frac{d}{4}\rceil}(\overline{\Omega};
\mathbb{R}^{d\times d})$ satisfies the ellipticity condition \eqref{ell}, $q \in C^{2\lceil\frac{d}{4}\rceil}(\overline{\Omega})$ is strictly positive on $\overline{\Omega}$.
\item[{\rm(ii)}] 
$f\in D(A^r)$ and $u_0\in D(A^{r+1})$, with $r>\frac{d+3}{2}$.
\end{itemize}
\end{definition}

We shall prove in Section 2 that for any admissible tuple $(\Omega,a,q,f,u_0)$ and $\sigma\in L^1(0,T)$, problem
\eqref{eqn:fde} has a unique weak solution  $u\in L^1(0,T; C^1(\overline{\Omega}))$. Further, we have
the following affirmative answers to the inverse problem for the cases $u_0\equiv0$ and $f\equiv0$,
respectively; for the detailed proofs, see Section \ref{sec:uniqueness}.
\begin{theorem}\label{thm:main1}
Let $(\Omega_k,a_k,q_k,f_k,0)$, $k=1,2$, be two admissible tuples with $u_0\equiv0$,
$\sigma\in L^1(0,T)$ be such that $\sigma\not\equiv0$, and the constants $r_1^k,\ldots,
r_{N_k}^k\in(0,+\infty)$, $0<\alpha_1^k<\ldots<\alpha_{N_k}^k<1$, $N_k\in\mathbb N$,
$k=1,2$. Let $u^k$ be the weak solution of problem \eqref{eqn:fde} with $(\Omega,a,q,f,
u_0)=(\Omega_k,a_k,q_k,f_k,0)$, $N=N_k$, $r_1=r_1^k$,$\ldots,$ $r_N=r_{N_k}^k$,
$\alpha_1=\alpha_1^k$,$\ldots,$ $\alpha_N=\alpha_{N_k}^k$.
Assume that there exist $x_k\in\partial\Omega_k$, $k=1,2$ such that
\begin{equation}\label{t1a}
\mathcal R_1^*f_1(x_1)\neq0,\quad \mathcal R_2^*f_2(x_2)\neq0,
\end{equation}
hold, and that one of the following conditions holds:
{\rm(i)} $\mathcal R_1^*f_1(x_1)=\mathcal R_2^*f_2(x_2)$ or {\rm(ii)} $r_{N_1}^1=r_{N_2}^2.$
Then the condition
\begin{equation}\label{t1b}
\mathcal R_1^*u^1(x_1,t)=\mathcal R_2^*u^2(x_2,t),\quad t\in(0,T)
\end{equation}
implies that $N_1=N_2=N$ and
\begin{equation}\label{t1c}
\mathcal R_1^*f_1(x_1)=\mathcal R_2^*f_2(x_2),\quad \alpha_1^1=\alpha_1^2,\ldots,\ \alpha_N^1=\alpha_N^2,\quad  r_1^1=r_1^2,\ldots,\ r_N^1=r_N^2.
\end{equation}
\end{theorem}

\begin{theorem}\label{thm:main2}
Let $(\Omega_k,a_k,q_k,0, u_0^k)$, $k=1,2$, be two admissible tuples with $f\equiv0$,  and the constants $r_1^k,\ldots,r_{N_k}^k\in(0,+\infty)$, $0<\alpha_1^k<\ldots<\alpha_{N_k}^k<1$, for $N_k\in\mathbb{N}$, $k=1,2$. Let $u^k$ be the weak solution of problem \eqref{eqn:fde} with $(\Omega,a,q,f,u_0)=(\Omega_k,a_k,q_k,0,u_0^k)$, $N=N_k$, $r_1=r_1^k$,$\ldots,$ $r_N=r_{N_k}^k$, $\alpha_1=\alpha_1^k$,
$\ldots,$ $\alpha_N=\alpha_{N_k}^k$. Assume that there exist  $x_k\in\partial\Omega_k$, $k=1,2$ such that the condition
\begin{equation}\label{t2a}
\mathcal R_1^*\mathcal A_1u^1_0(x_1)\neq0,\quad \mathcal R_2^*\mathcal A_2u^2_0(x_2)\neq0
\end{equation}
holds, and that one of the following conditions holds
{\rm(i)} $\mathcal R_1^*\mathcal A_1u_0^1(x_1)=\mathcal R_2^*\mathcal A_2u_0^2(x_2)$
or {\rm(ii)} $r_{N_1}^1=r_{N_2}^2$, where
$\mathcal A_k$, $k=1,2$, denotes the operator given by \eqref{A} with $a=a_k$, $q=q_k$ and $\Omega=\Omega_k$.
Then for any $T_1,T_2\in[0,T]$ satisfying $T_1<T_2$, the condition
\begin{equation}\label{ttt2a}
\mathcal R_1^*u^1(x_1,t)=\mathcal R_2^*u^2(x_2,t),\quad t\in(T_1,T_2)
\end{equation}
 implies that $N_1=N_2=N$ and conditions \eqref{t1c}  and
\begin{equation}\label{tt2a}
\mathcal R_1^*u_0^1(x_1)=\mathcal R_2^*u_0^2(x_2),\quad \mathcal R_1^*\mathcal A_1u_0^1(x_1)=\mathcal R_2^*\mathcal A_2u_0^2(x_2).
\end{equation}
\end{theorem}

Note that Theorems \ref{thm:main1} and \ref{thm:main2} are stated with unknown problem data in the sense that
they are stated with $\Omega_1\neq \Omega_2$, $a_1\neq a_2$, $q_1\neq q_2$, $f_1\neq f_2$ and $u_0^1\neq u_0^2$,
henceforth the term ``unknown medium''. Moreover, we mention that in these results the points $x_k\in\partial\Omega_k$,
$k=1,2$ do not need to coincide, and in fact we can even consider the case $\overline{\Omega_1}\cap \overline{\Omega_2}=
\emptyset$ or the case $\partial\Omega_1\cap \partial\Omega_2=\emptyset$. The main assumption that we impose on the data is
given by condition \eqref{t1a} or \eqref{t2a}, which require that the data $\mathcal R^* \mathcal Au_0$ and $\mathcal R^*f$ do
not vanish at the measurement point on $\partial\Omega$. Besides, we do not even need to assume that the
values at these points coincides provided that condition (ii) is fulfilled. The proof of the theorems relies on the solution
representation in the Laplace domain \cite[Section 6.2]{Jin:2021} and the time analyticity of the measurement. The uniqueness
proof is actually constructive, and motivates developing a simple recovery procedure based on the nonlinear least-squares, and asymptotic
expansion at $t=0$, cf. Proposition \ref{prop:asymptotic-u0}, directly inspired by the analysis. In Section \ref{sec:numer}, we present illustrative numerical
experiments to show the feasibility of the recovery using a nonlinear least-squares procedure (with either fractional
polynomials or fractional rational approximations as the regressor), when the measurement at small time is available.

Now we situate the uniqueness results in existing literature. The recovery of fractional orders probably
has been extensively studied; see \cite{LiLiuYamamoto:2019order}
for a survey. However, most existing studies focus on recovering
one single order in the model \eqref{eqn:fde0} \cite{AlimovAshurov:2020,ChengNakagawaYamamoto:2009,JinKian:2021,HNWY,LiaoWei:2019,
LiZhang:2020,Yamamoto:2021order}, sometimes together with other parameters, e.g.,
diffusion or potential coefficients,  given certain observational data. The only works on
recovering multiple orders are \cite{LiImanuvilovYamamoto:2016,LiYamamoto:2015,SunLiZhang:2021}. Li and
Yamamoto \cite{LiYamamoto:2015} proved the unique recovery of multiple orders in two cases: (i) the uniqueness in
simultaneously identifying $\{(\alpha_i,r_i)\}_{i=1}^N$ when $d=1$ and $u_0=\delta(x-x^*)$ (the Dirac delta function
concentrated at $x^*\in \Omega$) by measured data at one endpoint; (ii) the uniqueness in determining $\{(\alpha_i,
r_i)\}_{i=1}^N$ when $d\geq1$ and $u_0\in L^2(\Omega)$ by interior measurement. The analysis is based on the asymptotic
behavior of the multinomial Mittag-Leffler functions (cf. Remark \ref{rmk:asymptotic}). Li et al \cite{LiImanuvilovYamamoto:2016} proved the unique
recovery of orders and several coefficients from data consisting of a suitably defined Dirichlet-to-Neumann map.
Sun et al \cite{SunLiZhang:2021} proved the
unique recovery of the orders and the potential $q$ (for one-dimensional problem)
using the Gel'fand-Levitan theory for Sturm-Liouville problems. All these existing works assume
a fully known forward model. This work is also a natural continuation of the authors' recent
work \cite{JinKian:2021}, where the unique recovery of one single fractional order was proved for the
model \eqref{eqn:fde0} within an unknown medium (e.g. diffusion coefficient, potential) and scatterer
from lateral flux data at one single point (see also \cite{Yamamoto:2021order} in the case
of nonself-adjoint elliptic operators.) Note that the analysis \cite{JinKian:2021} relies
heavily on the analyticity of the solution at large time, asymptotic of the two-parameter
Mittag-Leffler function $E_{\alpha,\beta}(z)$ and the strong maximum principle and Hopf's lemma for
elliptic problems. Thus this work differs greatly from \cite{JinKian:2021} in the proof technique and
admissible observation data.

The rest of the paper is organized as follows. In Section \ref{sec:prelim} we recall preliminaries
of problem \eqref{eqn:fde}, e.g., existence, regularity and analyticity of the solution. Then in
Section \ref{sec:uniqueness}, we prove Theorems \ref{thm:main1} and \ref{thm:main2}.
Last, in Section \ref{sec:numer}, we present some numerical experiments for recovering the orders.
Throughout, the notation $\R_+$ denotes the set $(0,+\infty)$, By $\left\langle \cdot, \cdot\right\rangle$ we denote
the scalar product in $L^2(\Omega)$. The notation $C$ denotes a generic
constant which may change from one line to the next, but it is always independent of the quantity
under analysis, e.g. $p$.

\section{Preliminaries}\label{sec:prelim}

In this section we give several preliminary properties of problem \eqref{eqn:fde}, e.g., existence
of a weak solution and time analyticity. First we recall the concept of weak solutions.
Li et al \cite{LiLiuYamamoto:2015} proved the unique
existence of a mild solution using multinomial Mittag-Leffler functions. We employ
a representation of solutions in terms of inverse Laplace transform \cite{JinLiZhou:2018,KianSoccorsiXueYamamoto:2019,KianSoccorsiYamamoto:2018}
(or \cite[Section 6.2]{Jin:2021}).
\begin{definition}\label{def:weak}
Let $\sigma\in L^1(0,T)$ and $u_0,f\in L^2 (\Omega)$. A function $u$ is said to be a weak solution to
problem \eqref{eqn:fde} if there exists
$v\in L^1_{loc}(0,+\infty;L^2 (\Omega))$ satisfying
$u=v_{| \Omega\times(0,T)}$ and the following properties:
\begin{itemize}
\item[{\rm(i)}] $\inf\{\lambda>0:\ t\mapsto e^{-\lambda t}v(\cdot,t)\in
L^1(0,+\infty;L^2 (\Omega))\}=0$,
\item[{\rm(ii)}] for all $p>0$, the Laplace transform
$
\widehat{v}(\cdot, p):=\int_0^{+\infty}e^{-pt}v(\cdot,t)\d t
$
of $v$ belongs to $L^2 (\Omega)$ and solves
\begin{equation}\label{eqq2}
\left\{ \begin{aligned}
 \mathcal A \widehat v(p) +\sum_{k=1}^Nr_kp^{\alpha_k} \widehat v(p)& =  \Big(\int_0^T e^{-pt}\sigma(t)\d t\Big) f+\sum_{k=1}^Nr_kp^{\alpha_k-1}u_0,\quad
 \mbox{in  }\Omega,\\
\mathcal R \widehat v(p) & =  0, \quad \mbox{on } \partial \Omega.
\end{aligned}\right.
\end{equation}
\end{itemize}
\end{definition}

Throughout, we fix $\theta\in(\frac{\pi}{2},\pi)$, $\delta\in\R_+$ and the contour $\gamma(\delta,\theta)$ in $\mathbb C$ defined by
$\gamma(\delta,\theta):=\{\delta e^{i\beta}{\color{red}:}\ \beta\in[-\theta,\theta]\} \cup \gamma_\pm(\delta,\theta):=\{r e^{\pm i\theta}:\  r\geq\delta\}$
oriented in the counterclockwise direction. Let $\theta_1\in(0,\frac{2\theta-\pi}{8})$. For all $z\in D_{\theta_1}:=\{re^{i\beta}:\ r>0,\ \beta\in (-\theta_1,\theta_1)\}$,
we define two solution operators $S_1(z), S_2(z)\in\mathcal B(L^2(\Omega))$ by
\begin{align}
S_1(z)v&=\frac{1}{2i\pi}\int_{\gamma(\delta,\theta_1)} e^{z p}\Big(A+\sum_{k=1}^Nr_kp^{\alpha_k}\Big)^{-1}\Big(\sum_{k=1}^Nr_kp^{\alpha_k-1}\Big)v \d p,\quad v\in L^2(\Omega),\label{S1}\\
S_2(z)v&=\frac{1}{2i\pi}\int_{\gamma(\delta,\theta_1)}e^{z p}\Big(A+\sum_{k=1}^Nr_kp^{\alpha_k}\Big)^{-1}v \d p,\quad v\in L^2(\Omega).\label{S}
\end{align}
The operators $S_1(z)$ and $S_2(z)$ correspond to the initial data and the right
hand side, respectively.

Recall that the spectrum of the operator $A$
consists of a nondecreasing sequence of strictly positive eigenvalues
$(\lambda_{n})_{n\geq1}$ repeated with respect to their multiplicity.
In the Hilbert space $L^2(\Omega)$, we introduce an orthonormal basis
of eigenfunctions $(\phi_{n})_{n\geq1}$ of $A$ associated
with the eigenvalues $(\lambda_{n})_{n\geq1}$.
For all $s\geq 0$, we denote by $A^s$ the fractional power operator defined by
\[
A^s g=\sum_{n=1}^{+\infty}\langle g,\phi_n\rangle
\lambda_{n}^s\phi_{n},\quad g\in D(A^s)
= \Big\{{g}\in L^2(\Omega):\ \sum_{n=1}^{+\infty}\abs{\langle g,
\phi_{n}\rangle}^2 \lambda_{n}^{2s}<\infty
\Big\},
\]
and in $D(A^s)$, we define the graph norm $\|\cdot\|_{D(A^s)}$ by
\[\|g\|_{D(A^s)}
= \Big(\sum_{n=1}^{+\infty}\abs{\langle g,
\phi_{n}\rangle}^2 \lambda_{n}^{2s}\Big)^{\frac{1}{2}},
\quad g\in D(A^s).
\]
Following \cite[Lemma 3.4]{Kian:2020} and \cite[Theorem 1.2]{LiKianSoccorsi:2019}, we can prove the following result
\begin{lemma}\label{l1} For all $s\in[0,1]$, the map $z\mapsto S_j(z)$ is holomorphic in $D_{\theta_1}$ as a map taking values in $\mathcal B(L^2(\Omega); D(A^s))$ and there exists $C>0$ depending only on $\mathcal A$, $r_1,\ldots,r_N$, $\alpha_1,\ldots,\alpha_N$  and $\Omega$ such that
\begin{align*}
\norm{S_1(z)}_{B(L^2(\Omega); D(A^s))}&\leq C\max(|z|^{\alpha_1(1-s)-1},|z|^{\alpha_N(1-s)-1},1),\quad z\in D_{\theta_1},\\
\norm{S_2(z)}_{B(L^2(\Omega); D(A^s))}&\leq C\max(|z|^{-\alpha_1s},|z|^{-s\alpha_N},1),\quad z\in D_{\theta_1}.
\end{align*}
\end{lemma}

In a similar way to\cite[Proposition 2.1]{LiKianSoccorsi:2019}, one can prove that for
$\sigma\in L^\infty(0,T)$ and $u_0,f\in L^2(\Omega)$ problem \eqref{eqn:fde} admits a unique weak solution $u\in L^1(0,T;D(A^s))$, $s\in [0,1)$ given by
\begin{equation}\label{sol}
u(\cdot,t)=S_1(t)u_0+\int_0^t\sigma(s) S_2(t-s)f\d s,\quad t\in(0,T).
\end{equation}
We claim that the representation \eqref{sol} indeed gives a weak solution of \eqref{eqn:fde} in the sense of Definition \ref{def:weak}. We consider only the case $u_0\equiv0$. Since $\|(A+\sum_{k=1}^Nr_kp^{\alpha_k})^{-1}\|_{\mathcal B(L^2(\Omega))}\leq C(1+|p|)^{-\alpha_N}$,  $p\in\mathbb C\setminus(-\infty,0]$, we can define the operator-valued function
$$R_1(t):=  \frac{1}{2i\pi} \int_{-i\infty}^{+i\infty} e^{t p} (p+1)^{-1}\Big(A+\sum_{k=1}^Nr_k(p+1)^{\alpha_k}\Big)^{-1} \d p
, \quad t\in\R.$$
Following \cite[Proposition 2.1]{LiKianSoccorsi:2019}, $R_1\in L^\infty(\mathbb R; \mathcal B(L^2(\Omega)))$ is supported on $[0,+\infty)$. Moreover, using the argument of  \cite[Theorem 19.2 and the remark]{Ru}, we deduce that the Laplace transform $\widehat{R_1}(p)$ is well defined for $p>0$ and $\widehat{R_1}(p)=(p+1)^{-1}(A+\sum_{k=1}^Nr_k(p+1)^{\alpha_k})^{-1}$. Similarly, following \cite[Proposition 2.1]{LiKianSoccorsi:2019}, $R_2(t):=e^tR_1(t)$ belongs to the set of tempered distributions supported on $[0,+\infty)$ and taking values in $\mathcal B(L^2(\Omega))$. Let $\tilde{\sigma}$ be the extension of $\sigma$ to $\mathbb R$ by zero. Fixing the maps
$$w(\cdot,t)=(R_2(t))*(\tilde{\sigma}(t)f)=\int_0^t\tilde{\sigma}(s)R_2(t-s)f\d s,\quad v(\cdot,t)=\int_0^t\tilde{\sigma}(s) S_2(t-s)f\d s,\quad t\in\mathbb R$$
and repeating the arguments of \cite[Proposition 2.1]{LiKianSoccorsi:2019} give that $\partial_tw=v$ in the sense of tempered distributions. Therefore, $\widehat{v}(\cdot,p)=(A+\sum_{k=1}^Nr_kp^{\alpha_k})^{-1}\widehat{\tilde{\sigma}}(p)f$ is well defined for $p>0$ and it solves \eqref{eqq2}. Finally, for $u$ given by \eqref{sol}, we have $v=u$ on $\Omega\times(0,T)$, and hence the weak solution of \eqref{eqn:fde} takes the form \eqref{sol}.
These arguments show the existence of a weak solution of \eqref{eqn:fde} when $\sigma\in L^\infty(0,T)$. For the case $\sigma\in L^1(0,T)$, it suffices to combine Lemma \ref{l1} with a density argument (cf. \cite[Proposition 6.1]{KianSoccorsiXueYamamoto:2019}).
In a similar way to \cite{LiKianSoccorsi:2019}, one can check
$$S_j(z)h=\sum_{n=1}^{+\infty}S_{j,n}(z)\left\langle h,\phi_n\right\rangle\phi_n,\quad z\in D_{\theta_1},\ j=1,2,$$
with
\begin{align*}
S_{1,n}(z)&=\frac{1}{2i\pi}\int_{\gamma(\delta,\theta_1)} e^{z p}\Big(\lambda_n+\sum_{k=1}^Nr_kp^{\alpha_k}\Big)^{-1}\Big(\sum_{k=1}^Nr_kp^{\alpha_k-1}\Big) \d p,\quad z\in D_{\theta_1},\ n\in\mathbb N,\\
S_{2,n}(z)&=\frac{1}{2i\pi}\int_{\gamma(\delta,\theta_1)}e^{z p}\Big(\lambda_n+\sum_{k=1}^Nr_kp^{\alpha_k}\Big)^{-1} \d p,\quad z\in D_{\theta_1},\ n\in\mathbb N.
\end{align*}
In passing, note that the functions $S_{1,n}$ and $S_{2,n}$ can be expressed explicitly via
multinomial Mittag-Leffler functions, cf. Remark \ref{rmk:asymptotic}. Repeating the arguments
of Lemma \ref{l1}, we deduce that for all $n\in\mathbb N$ and $j=1,2$, $S_{j,n}$ is
holomorphic on $D_{\theta_1}$. Moreover, for all $s\in[0,1]$, we have
\begin{align}
\abs{S_{1,n}(z)}&\leq C\lambda_n^{-s}\max(|z|^{-s\alpha_1},|z|^{-s\alpha_N},1),\quad z\in D_{\theta_1},\ n\in\mathbb N,\label{l1d} \\
\abs{S_{2,n}(z)}&\leq C\lambda_n^{-s}\max(|z|^{\alpha_1(1-s)-1},|z|^{\alpha_N(1-s)-1},1),\quad z\in D_{\theta_1},\ n\in\mathbb N,\label{l1b}
\end{align}
with $C>0$ depending only on $\mathcal A$, $r_1,\ldots,r_N$, $\alpha_1,\ldots,\alpha_N$  and $\Omega$.
Using this result we can prove the following representation of the measured data.

\begin{lemma}\label{l2}
Let $f,u_0\in D(A^{\frac{d}{4}})$  and $\sigma\in L^1(0,T)$. Then the map $t\mapsto S_1(t)u_0$
and $t\mapsto S_2(t)f$  are analytic with respect to $t\in\R_+$ as a function taking values
in $C^1(\overline{\Omega})$. Moreover, problem \eqref{eqn:fde} admits a unique weak
solution $u\in L^1(0,T; C^1(\overline{\Omega}))$ satisfying
\begin{equation}\label{l2a}
\mathcal R^* u(x,t)=\mathcal R^* [S_1(t)u_0](x)+ \int_0^t\sigma(s) \mathcal R^* [S_2(t-s)f](x)\d s,\quad t\in(0,T),\ x\in\partial\Omega.
\end{equation}
\end{lemma}
\begin{proof}
Without loss of generality, we only prove the representation for $u_0\equiv0$.
In view of the identity \eqref{DA}, by interpolation, we deduce that the space $D(A^{\frac{d}{4}+\frac{3}{4}}) $
embeds continuously into $H^{\frac{d}{2}+\frac{3}{2}}(\Omega)$ and the Sobolev
embedding theorem implies that $D(A^{\frac{d}{4}+\frac{3}{4}}) $ embeds
continuously into $ C^1(\overline{\Omega})$. In addition, applying \eqref{l1b},
for all  $z\in D_{\theta_1}$ and all $m_1,m_2\in\mathbb N$, $m_1<m_2$, we have
\begin{equation}\label{l2b}
\begin{aligned}
\norm{\sum_{n=m_1}^{m_2}S_{2,n}(z)\left\langle f,\phi_n\right\rangle\phi_n}_{ C^1(\overline{\Omega})}&\leq C\norm{\sum_{n=m_1}^{m_2}S_{2,n}(z)\left\langle {f},\phi_n\right\rangle\phi_n}_{D(A^{\frac{d}{4}+\frac{3}{4}})}\\
\ &\leq C\max(|z|^{\frac{\alpha_1}{4}-1},|z|^{\frac{\alpha_N}{4}-1},1)\Big(\sum_{n=m_1}^{m_2}\lambda_n^{\frac{d}{2}}|\left\langle f,\phi_n\right\rangle|^2\Big)^{\frac{1}{2}},
\end{aligned}
\end{equation}
with $C>0$ independent of $z$, $m_1$ and $m_2$.
Combining this with the condition $f\in D(A^{\frac{d}{4}})$ yields that the sequence
$\sum_{n=1}^{N}S_{2,n}(z)\left\langle f,\phi_n\right\rangle\phi_n$, $N\in\mathbb N,$
converges uniformly with respect to $z$ on any compact set of $D_{\theta_1}$ to $S_2(z)f$
as a function taking values in $C^1(\overline{\Omega})$. This proves that the map
$D_{\theta_1}\ni z\mapsto S_2(z)f$ is holomorphic as a function taking values in
$C^1(\overline{\Omega})$. This implies the first statement of the lemma.
Next, in view of \eqref{l2b}, we have
$$
\norm{u(\cdot,t)}_{D(A^{\frac{d}{4}+\frac{3}{4}})}\leq C\norm{f}_{D(A^{\frac{d}{4}})}(\max(t^{\frac{\alpha_1}{4}-1},t^{\frac{\alpha_N}{4}-1},1)\mathds{1}_{(0,T)})*(|\sigma|\mathds{1}_{(0,T)})(t),\quad t\in(0,T),
$$
where $\mathds{1}_{(0,T)}$ denotes the characteristic function of $(0,T)$ and $*$ denotes
the convolution product. Therefore, applying Young's inequality, we obtain $u\in L^1(0,T;
D(A^{\frac{d}{4}+\frac{3}{4}}))\subset L^1(0,T; C^1(\overline{\Omega}))$, showing
the second assertion. In the
same way, applying \eqref{l2b}, we deduce \eqref{l2a}.
\end{proof}

\begin{lemma}\label{l3}
The following estimates hold
\begin{align}
\sum_{n=1}^\infty\lambda_n|\left\langle v,\phi_n\right\rangle|\norm{\mathcal R^*\phi_n}_{L^\infty(\partial\Omega)}\leq C\norm{v}_{D(A^s)},\quad \forall v\in D(A^s), s>\frac{d}{2}+\frac{3}{2}, \label{l3a}\\
\mathcal R^*v(x)=\sum_{n=1}^\infty \left\langle v,\phi_n\right\rangle \mathcal R^*\phi_n(x),\quad  x\in\partial\Omega,\quad \forall v\in D(A^s), s>\frac{d}{2}+\frac{1}{2}.\label{l3b}
\end{align}
\end{lemma}
\begin{proof}
Observe that, according to the Weyl's asymptotic
formula \cite{Weyl:1912}, there exists $C>0$ such that
$C^{-1} n^{\frac{2}{d}}\leq\lambda_n\leq C n^{\frac{2}{d}}$, for all $n\geq1.$
Thus, we obtain for any $r>\frac d2$
\begin{equation}\label{l3c}
\sum_{n=1}^\infty \lambda_n^{-r}\leq C\sum_{n=1}^\infty n^{-\frac{2}{d}r}<\infty.
\end{equation}
Meanwhile, the Sobolev embedding theorem implies for any $n\in\mathbb{N}$ and $\epsilon>0$,
\begin{align*}
\norm{\mathcal R^*\phi_n}_{L^\infty(\partial\Omega)}&\leq
   C\norm{\phi_n}_{C^1(\overline{\Omega})} {\leq C\norm{\phi_n}_{D(A^{\frac{d}{4}+\frac12+\epsilon})}\leq C\lambda_n^{\frac{d}{4}+\frac{1}{2}+\epsilon},}
\end{align*}
Then the Cauchy-Schwarz inequality implies
\begin{align*}
&\sum_{n=1}^\infty\lambda_n|\langle v,\phi_n\rangle|\norm{\mathcal R^*\phi_n}_{L^\infty(\partial\Omega)}\leq C
\sum_{n=1}^\infty\lambda_n^{\frac{d}{4}+\frac{3}{2}+\epsilon}|\left\langle v,\phi_n\right\rangle|\\
= &C\sum_{n=1}^\infty\lambda_n^{\frac{d}{4}+\frac{3}{2}+\epsilon+\frac r2}|\langle v,\phi_n\rangle|\lambda_n^{-\frac r2}
\leq C\Big(\sum_{n=1}^\infty\lambda_n^{\frac d2+3+2\epsilon+r}|\langle v,\phi_n\rangle|^{2}\Big)^{\frac{1}{2}}\Big(\sum_{n=1}^\infty\lambda_n^{-r}\Big)^{\frac{1}{2}}.
\end{align*}
Combining this with \eqref{l3c} and the condition $v \in D(A^s)$ gives \eqref{l3a}. This argument
also shows that $\sum_{n=1}^{N}\left\langle v,\phi_n\right\rangle\phi_n$, $N\in\mathbb N,$
converges in $C^1(\overline{\Omega})$ to $v$, which directly gives \eqref{l3b}.
\end{proof}

\begin{remark}
The regularity on $v$ in \eqref{l3a} can be relaxed to $v\in A^s$, $s>\frac{d}{2}+1$,
if $\mathcal{R}^*$ is the Dirichlet trace operator, and a similar observation holds
for the estimate \eqref{l3b}.
\end{remark}

\section{Proof of Theorems \ref{thm:main1} and \ref{thm:main2}}\label{sec:uniqueness}

In this section, we give the proof of Theorems \ref{thm:main1} and \ref{thm:main2}.

\subsection{Proof of Theorem \ref{thm:main1}}
Throughout this part, the assumption of Theorem \ref{thm:main1} is fulfilled and prove that \eqref{t1b} implies that $N_1=N_2$ and \eqref{t1c} is fulfilled. Let $S_2^j$ correspond to \eqref{S} with $N=N_j$, $r_k=r^j_k$, $\alpha_k=\alpha^j_k$, $k=1,\ldots,N_j$
and with $A=A_j=\mathcal A_j$, acting in $L^2(\Omega_j)$ and with the boundary condition given by $\mathcal
R=\mathcal R_j$. We divide the lengthy proof into three steps.\\
\noindent \textbf{Step 1.} In this step, we show that \eqref{t1b} implies
\begin{equation}\label{t1d}
\mathcal R_1^*\Big[\Big(A_1+\sum_{k=1}^{N_1}r_k^1p^{\alpha_k^1}\Big)^{-1}f_1\Big](x_1)=\mathcal R_2^*\Big[\Big(A_2+\sum_{k=1}^{N_2}r_k^2p^{\alpha_k^2}\Big)^{-1}f_2\Big](x_2),\quad p\in\R_+.
\end{equation}
By Lemma \ref{l2}, we have
$$
\mathcal R^*_j u_j(x_j,t)=\int_0^t\sigma(s) \mathcal R^*_j [S_2^j(t-s)f_j](x_j)\d s,\quad t\in(0,T),\ j=1,2.
$$
 Let $v_j(t)= \mathcal R^*_j [S_2^j(t)f_j](x_j)$, $j=1,2$,
as a function in $L^1(0,T)$, cf. Lemma \ref{l2}. Therefore, condition \eqref{t1b} implies
$\int_0^t\sigma(s)[v_1(t-s)-v_2(t-s)]\d s=0$, for $t\in(0,T).$
By Titchmarsh convolution theorem \cite[Theorem VII]{Titchmarsh:1926},
there exist $T_1,T_2\in[0,T]$ such that $T_1+T_2\geq T$, $\sigma_{|(0,T_1)}\equiv0$ and $(v_1-v_2)_{|(0,T_2)}
\equiv0$. Meanwhile, since $\sigma\not\equiv0$, $T_1<T$. Thus, $T_2=T_1+T_2-T_1\geq T-T_1>0$ and
$$
\mathcal R^*_1 [S_2^1(t)f_1](x_1)=\mathcal R^*_2 [S_2^2(t)f_2](x_2),\quad t\in(0,T_2).
$$
The analyticity of the maps $\R_+\ni t\mapsto \mathcal R^*_j [S_2^j(t)f_j](x_j)$, $j=1,2$, given in Lemma \ref{l2}, implies
\begin{equation}\label{t1e}
\mathcal R^*_1 [S_2^1(t)f_1](x_1)=\mathcal R^*_2 [S_2^2(t)f_2](x_2),\quad t\in\R_+.
\end{equation}
Moreover, applying the properties of the map \eqref{S} given in Section \ref{sec:prelim}, we have
$$
\widehat{S_2^j(t)f_j}(p)=\Big(A_j+\sum_{k=1}^{N_j}r_k^jp^{\alpha_k^j}\Big)^{-1}f_j,
$$
and by the arguments of Lemma \ref{l2}, we obtain
$$
\widehat{\mathcal R^*_jS_2^j(t)f_j(\cdot,x_j)}(p)=\mathcal R^*_j\Big[\Big(A_j+\sum_{k=1}^{N_j}r_k^jp^{\alpha_k^j}\Big)^{-1}f_j\Big](x_j).
$$
Combining this with \eqref{t1e} leads to \eqref{t1d}.

\noindent\textbf{Step 2.} Now we fix $N=\min(N_1,N_2)$ and prove that condition \eqref{t1d} implies
\begin{equation}\label{t1f}
\left\{\begin{aligned}
 & \mathcal R_1^*f_1(x_1)=\mathcal R_2^*f_2(x_2),\\
 &\alpha_{N_1}^1=\alpha_{N_2}^2,\ldots,\alpha_{N_1-N+1}^1=\alpha_{N_2-N+1}^2, \\
 &r_{N_1}^1=r_{N_2}^2,\ldots,r_{N_1-N+1}^1=r_{N_2-N+1}^2.
\end{aligned}\right.
\end{equation}
We prove this result iteratively using the asymptotic properties of
$$
\mathcal R_j^*\Big(A_j+\sum_{k=1}^{N_j}r_k^jp^{\alpha_k^j}\Big)^{-1}f_j(x_1),\quad j=1,2,\ p\to+\infty.
$$
To this end, we fix $A=A_j$ with $\mathcal A=\mathcal A_j$, $\mathcal R=\mathcal R_j$, $\Omega=\Omega_j$,
and denote the non-decreasing sequence of strictly positive eigenvalues
of the operator $A_j$ by $(\lambda_{n}^j)_{n\geq1}$ and an $L^2(\Omega_j)$ orthonormal basis
of eigenfunctions $(\phi_{n}^j)_{n\geq1}$ associated
with the eigenvalues $(\lambda_{n}^j)_{n\geq1}$. Then, we have
$$
\Big(A_j+\sum_{k=1}^{N_j}r_k^jp^{\alpha_k^j}\Big)^{-1}f_j=\sum_{n=1}^\infty\frac{\langle f_j,\phi_n^j\rangle_{L^2(\Omega_j)} }{\lambda_n^j+\sum_{k=1}^{N_j}r_k^jp^{\alpha_k^j}}\phi_{n}^j,\quad p\in\R_+, j=1,2.
$$
Since {$f_j\in D(A_j^r)$, $r>\frac{d+3}{2}$, we have $(A_j+\sum_{k=1}^{N_j}r_k^jp^{\alpha_k^j})^{-1}f_j\in D(A_j^r)$.}
Moreover, following the argument of Lemma \ref{l3}, we obtain
$$\mathcal R_j^*\Big(A_j+\sum_{k=1}^{N_1}r_k^jp^{\alpha_k^j}\Big)^{-1}f_j(x_j)=\sum_{n=1}^\infty\frac{\langle f_j,\phi_n^j\rangle_{L^2(\Omega_j)} \mathcal R_j^*\phi_{n}^j(x_j)}{\lambda_n^j+\sum_{k=1}^{N_j}r_k^jp^{\alpha_k^j}},\quad p\in\R_+, j=1,2.$$
Combining this with \eqref{t1d} gives
\begin{equation}\label{t1g}
\sum_{n=1}^\infty\frac{\langle f_1,\phi_n^1\rangle_{L^2(\Omega_1)} \mathcal R_1^*\phi_{n}^1(x_1)}{\lambda_n^1+\sum_{k=1}^{N_1}r_k^1p^{\alpha_k^1}}=     \sum_{n=1}^\infty\frac{\langle f_2,\phi_n^2\rangle_{L^2(\Omega_2)} \mathcal R_2^*\phi_{n}^2(x_2)}{\lambda_n^2+\sum_{k=1}^{N_2}r_k^2p^{\alpha_k^2}},\quad p\in\R_+.
\end{equation}
Meanwhile, since $f_j\in D(A_j^r)$, $j=1,2$, by the mean value theorem, we deduce
\begin{align*}
&\sum_{n=1}^\infty\frac{\langle f_j,\phi_n^j\rangle_{L^2(\Omega_j)} \mathcal R_j^*\phi_{n}^j(x_j)}{\lambda_n^j+\sum_{k=1}^{N_j}r_k^jp^{\alpha_k^j}}\\
&=\sum_{n=1}^\infty\frac{\langle f_j,\phi_n^j\rangle_{L^2(\Omega_j)} \mathcal R_j^*\phi_{n}^j(x_j)}{\sum_{k=1}^{N_j}r_k^jp^{\alpha_k^j}}-\sum_{n=1}^\infty\int_0^1\frac{\lambda_n^j\langle f_j,\phi_n^j\rangle_{L^2(\Omega_j)} \mathcal R_j^*\phi_{n}^j(x_j)}{(\sum_{k=1}^{N_j}r_k^jp^{\alpha_k^j}+s\lambda_n^j)^2}\d s.
\end{align*}
By Lemma \ref{l3}, there holds
$$
\mathcal R_j^*f_j(x_j)=\sum_{n=1}^\infty\langle f_j,\phi_n^j\rangle_{L^2(\Omega_j)} \mathcal R_j^*\phi_{n}^j(x_j),
$$
and thus we obtain
$$
\sum_{n=1}^\infty\frac{\langle f_j,\phi_n^j\rangle_{L^2(\Omega_j)} \mathcal R_j^*\phi_{n}^j(x_j)}{\lambda_n^j+\sum_{k=1}^{N_j}r_k^jp^{\alpha_k^j}}=\frac{ \mathcal R_j^*f_j(x_j)}{\sum_{k=1}^{N_j}r_k^jp^{\alpha_k^j}}-\sum_{n=1}^\infty\int_0^1\frac{\lambda_n^j\langle f_j,\phi_n^j\rangle_{L^2(\Omega_j)} \mathcal R_j^*\phi_{n}^j(x_j)}{(\sum_{k=1}^{N_j}r_k^jp^{\alpha_k^j}+s\lambda_n^j)^2}\d s.
$$
Then, it follows that, for $p>1$, there holds
$$\abs{\sum_{n=1}^\infty\frac{\langle f_j,\phi_n^j\rangle_{L^2(\Omega_j)} \mathcal R_j^*\phi_{n}^j(x_j)}{\lambda_n^j+\sum_{k=1}^{N_j}r_k^jp^{\alpha_k^j}}-\frac{ \mathcal R_j^*f_j(x_j)}{\sum_{k=1}^{N_j}r_k^jp^{\alpha_k^j}}}\leq C\frac{\sum_{n=1}^\infty\lambda_n^j|\langle f_j,\phi_n^j\rangle_{L^2(\Omega_j)} \mathcal R_j^*\phi_{n}^j(x_j)|}{(\sum_{k=1}^{N_j}r_k^jp^{\alpha_k^j})^2}.$$
Thus, Lemma \ref{l3} implies
$$\abs{\sum_{n=1}^\infty\frac{\langle f_j,\phi_n^j\rangle_{L^2(\Omega_j)} \mathcal R_j^*\phi_{n}^j(x_j)}{\lambda_n^j+\sum_{k=1}^{N_j}r_k^jp^{\alpha_k^j}}-\frac{ \mathcal R_j^*f_j(x_j)}{\sum_{k=1}^{N_j}r_k^jp^{\alpha_k^j}}}\leq C\frac{\|f_j\|_{D(A_j^r)}}{(\sum_{k=1}^{N_j}r_k^jp^{\alpha_k^j})^2}.$$
Therefore, for $j=1,2$, we have
$$\sum_{n=1}^\infty\frac{\langle f_j,\phi_n^j\rangle_{L^2(\Omega_j)} \mathcal R_j^*\phi_{n}^j(x_j)}{\lambda_n^j+\sum_{k=1}^{N_j}r_k^jp^{\alpha_k^j}}=\frac{ \mathcal R_j^*f_j(x_j)}{\sum_{k=1}^{N_j}r_k^jp^{\alpha_k^j}}+\underset{p\to+\infty}{\mathcal O}\Big(\frac{1}{(\sum_{k=1}^{N_j}r_k^jp^{\alpha_k^j})^2}\Big)$$
and by combining this with \eqref{t1g}, we obtain
\begin{equation}\label{tt1g}
\frac{ \mathcal R_1^*f_1(x_1)}{\sum_{k=1}^{N_1}r_k^1p^{\alpha_k^1}}+\underset{p\to+\infty}{\mathcal O}\left(\frac{1}{(\sum_{k=1}^{N_1}r_k^1p^{\alpha_k^1})^2}\right)=    \frac{ \mathcal R_2^*f_2(x_2)}{\sum_{k=1}^{N_2}r_k^2p^{\alpha_k^2}}+\underset{p\to+\infty}{\mathcal O}\left(\frac{1}{(\sum_{k=1}^{N_2}r_k^2p^{\alpha_k^2})^2}\right).
\end{equation}
Next we use this identity to prove \eqref{t1f}. We start by proving
\begin{equation}\label{t1h}
\alpha_{N_1}^1=\alpha_{N_2}^2,\quad  r_{N_1}^1=r_{N_2}^2,\quad \mathcal R_1^*f_1(x_1)=\mathcal R_2^*f_2(x_2).
\end{equation}
Recall that, for $j=1,2$,
$$
\frac{ \mathcal R_j^*f_j(x_j)}{\sum_{k=1}^{N_j}r_k^jp^{\alpha_k^j}}=\frac{ \mathcal R_j^*f_j(x_j)}{r_{N_j}^jp^{\alpha_{N_j}^j}}+\underset{p\to+\infty}{\mathcal O}\left(p^{\alpha_{N_j-1}^j-2\alpha_{N_j}^j}\right).
$$
Then we deduce from \eqref{tt1g} that
$$
\frac{ \mathcal R_1^*f_1(x_1)}{r_{N_1}^1p^{\alpha_{N_1}^1}}+\underset{p\to+\infty}{\mathcal O}\left(p^{\alpha_{N_1-1}^1-2\alpha_{N_1}^1}\right)=    \frac{ \mathcal R_2^*f_2(x_2)}{r_{N_2}^2p^{\alpha_{N_2}^2}}+\underset{p\to+\infty}{\mathcal O}\left(p^{\alpha_{N_2-1}^1-2\alpha_{N_2}^2}\right).
$$
In view of this identity, the fact that \eqref{t1a} is fulfilled and the fact that $\alpha_{N_j-1}^j<\alpha_{N_j}^j$, we deduce
$$
\alpha_{N_1}^1=\alpha_{N_2}^2,\quad \mathcal R_1^*f_1(x_1)(r_{N_1}^1)^{-1}= \mathcal R_2^*f_2(x_2)(r_{N_2}^2)^{-1}.
$$
Combining this with  condition (i) or (ii) of Theorem \ref{thm:main1}, we obtain \eqref{t1h}.
Now assume that there exists $\ell\in\{0,\ldots,N-2\}$ such that the following condition is fulfilled
\begin{equation}\label{t1j}
\mathcal R_1^*f_1(x_1)=\mathcal R_2^*f_2(x_2),\ \alpha_{N_1}^1=\alpha_{N_2}^2,\ldots,\alpha_{N_1-\ell}^1=\alpha_{N_2-\ell}^2,\  r_{N_1}^1=r_{N_2}^2,\ldots,r_{N_1-\ell}^1=r_{N_2-\ell}^2.
\end{equation}
We claim that this condition implies
\begin{equation}\label{t1k}
\alpha_{N_1}^1=\alpha_{N_2}^2,\ldots,\alpha_{N_1-\ell-1}^1=\alpha_{N_2-\ell-1}^2,\quad r_{N_1}^1=r_{N_2}^2,\ldots,r_{N_1-\ell-1}^1=r_{N_2-\ell-1}^2.
\end{equation}
Indeed, by letting
$
Q_\ell^j(p)=\sum_{k=0}^{\ell}r_{N_j-k}^jp^{\alpha_{N_j-k}^j},
$
for $j=1,2$, we find
$$
\frac{ \mathcal R_j^*f_j(x_j)}{\sum_{k=1}^{N_j}r_k^jp^{\alpha_k^j}}=\frac{ \mathcal R_j^*f_j(x_j)}{Q_\ell^j(p)}-\frac{\mathcal R_j^*f_j(x_j)r_{N_j-\ell-1}^jp^{\alpha_{N_j-\ell-1}^j}}{Q_\ell^j(p)^2}+\underset{p\to+\infty}{\mathcal O}\left(\frac{p^{\alpha_{N_j-\ell-2}^j}}{Q_\ell^j(p)^2}\right),
$$
where $\alpha^j_0=0$. Combining this with \eqref{tt1g} leads to
\begin{align*}
&\frac{ \mathcal R_1^*f_1(x_1)}{Q_\ell^1(p)}-\frac{\mathcal R_1^*f_1(x_1)r_{N_1-\ell-1}^1p^{\alpha_{N_1-\ell-1}^1}}{Q_\ell^1(p)^2}+\underset{p\to+\infty}{\mathcal O}\left(\frac{p^{\alpha_{N_1-\ell-2}^1}}{Q_\ell^1(p)^2}\right)\\
=&\frac{ \mathcal R_2^*f_2(x_2)}{Q_\ell^2(p)}-\frac{\mathcal R_2^*f_2(x_2)r_{N_2-\ell-1}^2p^{\alpha_{N_2-\ell-1}^2}}{Q_\ell^2(p)^2}+\underset{p\to+\infty}{\mathcal O}\left(\frac{p^{\alpha_{N_2-\ell-2}^2}}{Q_\ell^2(p)^2}\right).
\end{align*}
Further, condition \eqref{t1j} implies that, for all $p>0$, $Q_\ell^1(p)=Q_\ell^2(p)$ and $\mathcal R_1^*f_1(x_1)=\mathcal R_2^*f_2(x_2)$. Therefore, this identity and conditions \eqref{t1a} and \eqref{t1d} imply the claim \eqref{t1k}. Combining the iteration argument with \eqref{t1h} implies that  \eqref{t1f} holds.

\noindent\textbf{Step 3.} In this step we assume that \eqref{t1f} holds and complete the proof by proving $N_1=N_2$. Indeed, assuming that $N_1\neq N_2$, we may assume $N_1<N_2$. Then, using \eqref{t1a}, \eqref{t1f} and fixing
 $$Q_{N_1}(p)=\sum_{k=0}^{N_1-1}r_{N_1-k}^1p^{\alpha_{N_1-k}^1}=\sum_{k=0}^{N_1-1}r_{N_2-k}^2p^{\alpha_{N_2-k}^2},\quad b=\mathcal R_1^*f_1(x_2)=\mathcal R_2^*f_2(x_2)\neq0,$$
and,  repeating the arguments of Step 2, we deduce that
$$\frac{b}{Q_{N_1}(p)}-\frac{br_{N_2-N_1}^2p^{\alpha_{N_2-N_1}^2}}{Q_{N_1}(p)^2}+\underset{p\to+\infty}{\mathcal O}\left(\frac{p^{\alpha_{N_2-N_1-1}^2}}{Q_{N_1}(p)^2}\right)=\frac{b}{Q_{N_1}(p)}+\underset{p\to+\infty}{\mathcal O}\left(\frac{1}{Q_{N_1}(p)^2}\right),$$
with the convention $\alpha^2_0=0$.
Since $b\neq0$, $r_{N_2-N_1}^2>0$  and $\alpha_{N_2-N_1}^2>0$, the above identity cannot be true. This leads to a contradiction and $N_1=N_2$. Therefore, condition \eqref{t1f} implies \eqref{t1c}. This completes the proof of Theorem \ref{thm:main1}.

\subsection{Proof of Theorem \ref{thm:main2}}
Assume that  there exist $T_1,T_2\in[0,T]$, with $T_1<T_2$, such that \eqref{ttt2a} is fulfilled and we show that
$N_1=N_2$ and conditions \eqref{t1c} and \eqref{tt2a} are fulfilled under the assumptions of Theorem \ref{thm:main1}.
Following the argumentation of Section 2.1, we deduce that, for $j=1,2$, we have $u_j(\cdot,t)=S_1^j(t)u_0$, with
$S_1^j$ corresponding to \eqref{S1}  with $N=N_j$, $r_k=r^j_k$, $\alpha_k=\alpha^j_k$, $k=1,\ldots,N_j$ and with
$A=A_j=\mathcal A_j$, acting in $L^2(\Omega_j)$ and with the boundary condition given by $\mathcal R=\mathcal R_j$. Then, condition \eqref{ttt2a} implies
\begin{equation*}
\mathcal{R}_1^*[S_1^1(t)u^1_0](x_1,t)=\mathcal{R}_2^*[S_1^2(t)u^2_0](x_2,t),\quad t\in(T_1,T_2).
\end{equation*}
This and the analiticity of the maps $\R_+\ni t\mapsto \mathcal R^*_j [S_1^j(t)u^j_0](x_j)$,
$j=1,2$, cf. Lemma \ref{l2}, give
\begin{equation*}
\mathcal{R}_1^*[S_1^1(t)u^1_0](x_1,t)=\mathcal{R}_2^*[S_1^2(t)u^2_0](x_2,t),\quad t\in\R_+.
\end{equation*}
Therefore, repeating the argument in the first step of Theorem \ref{thm:main1}, we deduce that \eqref{ttt2a} implies
$$\begin{aligned}
&p^{-1}\mathcal R_1^*\Big[\Big(A_1+\sum_{k=1}^{N_1}r_k^1p^{\alpha_k^1}\Big)^{-1}\Big(\sum_{k=1}^{N_1}r_k^1p^{\alpha_k^1}\Big)u_0^1\Big](x_1)\\
=&\mathcal R_1^*\Big[\Big(A_1+\sum_{k=1}^{N_1}r_k^1p^{\alpha_k^1}\Big)^{-1}\Big(\sum_{k=1}^{N_1}r_k^1p^{\alpha_k^1-1}\Big)u_0^1\Big](x_1)\\
=&\mathcal R_2^*\Big[\Big(A_2+\sum_{k=1}^{N_2}r_k^2p^{\alpha_k^2}\Big)^{-1}\Big(\sum_{k=1}^{N_1}r_k^2p^{\alpha_k^2-1}\Big)u_0^2\Big](x_2)\\
=&p^{-1}\mathcal R_2^*\Big[\Big(A_2+\sum_{k=1}^{N_2}r_k^2p^{\alpha_k^2}\Big)^{-1}\Big(\sum_{k=1}^{N_1}r_k^2p^{\alpha_k^2}\Big)u_0^2\Big](x_2).\end{aligned}$$
Multiplying both sides of this identity by $p$, we obtain
\begin{equation}\label{t2d}
\begin{aligned}
&\mathcal R_1^*u_0^1(x_1)-\mathcal R_1^*\Big[\Big(A_1+\sum_{k=1}^{N_1}r_k^1p^{\alpha_k^1}\Big)^{-1}A_1u_0^1\Big](x_1)\\
=&\mathcal R_2^*u_0^2(x_2)-\mathcal R_2^*\Big[\Big(A_2+\sum_{k=1}^{N_2}r_k^2p^{\alpha_k^2}\Big)^{-1}A_2u_0^2\Big](x_2),\quad p\in\R_+.
\end{aligned}
\end{equation}
Moreover, repeating the arguments of the preceding section we deduce
$$
\lim_{p\to+\infty}\mathcal R_j^*\Big[\Big(A_j+\sum_{k=1}^{N_j}r_k^jp^{\alpha_k^j}\Big)^{-1}A_ju_0^j\Big](x_j)=0,\quad j=1,2.
$$
Therefore, condition \eqref{t2d} implies that $\mathcal R_1^*u_0^1(x_1)=\mathcal R_2^*u_0^2(x_2)$ and consequently,
\begin{equation}\label{t2e}
\mathcal R_1^*\Big[\Big(A_1+\sum_{k=1}^{N_1}r_k^1p^{\alpha_k^1}\Big)^{-1}A_1u_0^1\Big](x_1)=\mathcal R_2^*\Big[\Big(A_2+\sum_{k=1}^{N_2}r_k^2p^{\alpha_k^2}\Big)^{-1}A_2u_0^2\Big](x_2),\quad p\in\R_+.
\end{equation}
Combining this identity with the arguments of Step 2 of Theorem \ref{thm:main1} leads to
\begin{equation}\label{id}
\frac{\mathcal R_1^*A_1u_0^1(x_1)}{\sum_{k=1}^{N_1}r_k^1p^{\alpha_k^1}}+\underset{p\to+\infty}{\mathcal O}\Big(\frac{1}{(\sum_{k=1}^{N_1}r_k^1p^{\alpha_k^1})^2}\Big)=    \frac{ \mathcal R_2^*A_2u_0^2(x_2)}{\sum_{k=1}^{N_2}r_k^2p^{\alpha_k^2}}+\underset{p\to+\infty}{\mathcal O}\Big(\frac{1}{(\sum_{k=1}^{N_2}r_k^2p^{\alpha_k^2})^2}\Big).
\end{equation}
Therefore, repeating the arguments used in Steps 2 and 3 of Theorem \ref{thm:main1}, we deduce that
$N_1=N_2$, and conditions \eqref{t1c} and \eqref{tt2a} hold.

\begin{remark}
In view of the $t$-analyticity of the solution $u(x_0,t)$, the assertion in Theorem \ref{thm:main2} remains
valid, if the time trace $\mathcal{R}^*u(x_0,t)$ is only observed at a countable discrete set $t_1<t_2<\ldots$
with an accumulation point in the open interval $(T_1,T_2)$.
\end{remark}

\begin{remark}
Theorem \ref{thm:main2} can also be seen as follows. The solution $u$ of problem
\eqref{eqn:fde} with $\sigma\equiv0$ takes the form $u=u_0+v$, with $v$ solving
$$\left\{\begin{aligned}
\sum_{j=1}^Nr_j\partial_t^{\alpha_j}v +\mathcal{A} v &=  -\mathcal{A} u_0, \quad \mbox{in }\Omega\times(0,+\infty),\\
\mathcal R v&= 0, \quad \mbox{on } \partial\Omega\times(0,+\infty), \\
v&=0, \quad \mbox{in } \Omega\times \{0\}.
\end{aligned}\right.$$
In terms of Laplace transform, this implies $$\hat{u}(p)=\frac{u_0}{p}+\hat{v}(p)=p^{-1}\Big[u_0-\Big(A+\sum_{k=1}^{N_1}r_kp^{\alpha_k}\Big)^{-1}Au_0\Big],$$
which leads directly to formula \eqref{t2d}.
\end{remark}

\section{Numerical experiments and discussions}\label{sec:numer}

In this section, we illustrate the feasibility of the recovery of orders and the associated weights
with a set of numerical experiments, and discuss the potential pitfalls.

\subsection{Asymptotic expansion}
First we develop an algorithm for the numerical recovery, inspired by the analysis. Note that to recover the orders
$\alpha_j$ and weights $r_j$, one classical approach is to apply the regularization method,
e.g., Tikhonov regularization \cite{ItoJin:2015}, which involves a misfit on the measured
data (and proper regularization), with the forward map defined implicitly by problem
\eqref{eqn:fde} (as done recently in \cite{LiaoWei:2019,SunLiZhang:2021} for recovering
one single order). Unfortunately, this approach does not apply in the setting of this work,
since the medium (and thus the forward map) is unknown. Instead, we  employ a
more direct approach, which is inspired by the uniqueness analysis in Section \ref{sec:uniqueness}.
In the spirit of the classical Karamata-Feller Tauberian theorem \cite[Section XIII.5]{Feller:1970}, the
asymptotic of the Laplace transform $\hat u(x_0,p)$ as $p\to\infty$ corresponds to the asymptotic
of the function $u(x_0,t)$ as $t\to0^+$. Indeed, the argument for the uniqueness result in Theorem
\ref{thm:main2} is essentially about the asymptotic behavior of the measured data $\mathcal{R}^*
u(x_0,t)$. We have the following asymptotic expansion, which lays the foundation
of the procedure for the numerical recovery.

\begin{proposition}\label{prop:asymptotic-u0}
Let $(\Omega,a,q,f,u_0)$ be an admissible tuple. If $u_0\not\equiv0$ and $f\equiv0$, the solution $u$ to problem \eqref{eqn:fde} satisfies the
following asymptotic:
\begin{equation*}
  \mathcal{R}^* u(x_0,t)
                = \mathcal{R}^* u_0(x_0) - \mathcal{R}^*Au_0(x_0)\Big(r_N\frac{t^{\alpha_N}}{\Gamma(\alpha_N+1)}-\sum_{i=1}^{N-1}\frac{r_it^{2\alpha_N-\alpha_i}}{r_N^2\Gamma(2\alpha_N-\alpha_i+1)}\Big) + \mathcal{O}(t^{2\alpha_N}).
\end{equation*}
Similarly, if $u_0\equiv0$ and $f\not\equiv0$ and $\sigma\in L^1(0,T)$ with
$$\hat\sigma(p):=\mathcal{L}[\sigma](p)=c_0p^{-a-1} + \mathcal{O}(p^{-a-2}) \quad \mbox{as }p\to\infty,$$
for some $a\in[0,1]$, then there holds
\begin{equation*}
  \mathcal{R}^* u(x_0,t)= c_0\mathcal{R}^*f(x_0)\Big(r_N\frac{t^{\alpha_N+a}}{\Gamma(\alpha_N+a+1)}-\sum_{i=1}^{N-1}\frac{r_it^{2\alpha_N-\alpha_i+a}}{r_N^2\Gamma(2\alpha_N-\alpha_i+a+1)}\Big) + \mathcal{O}(t^{2\alpha_N}).
\end{equation*}
\end{proposition}
\begin{proof}
Indeed, the proof of Theorem \ref{thm:main2} gives the following
expansion for $u(x_0,p)$:
\begin{align*}
  &\mathcal{R}^*\hat u(x_0,p)  = \sum_{n=1}^\infty \frac{\sum_{i=1}^Nr_ip^{\alpha_i-1}}{\lambda_n+\sum_{i=1}^{N}r_ip^{\alpha_i}}\langle u_0,\varphi_n\rangle\mathcal{R}^*\varphi_n(x_0)\\
                 =& p^{-1}\sum_{n=1}^\infty\langle u_0,\varphi_n\rangle\mathcal{R}^*\varphi_n(x_0) -  \sum_{n=1}^\infty\frac{\lambda_np^{-1}}{\lambda_n+\sum_{i=1}^{N}r_ip^{\alpha_i}}\langle u_0,\varphi_n\rangle\mathcal{R}^*\varphi_n(x_0).
\end{align*}
In view of the following identity for large $p$
\begin{align*}
  &\frac{1}{\lambda + \sum_{i=1}^Nr_ip^{\alpha_i}} = \frac{1}{r_Np^{\alpha_N}}\frac{1}{1+\frac{\lambda+\sum_{i=1}^{N-1}r_ip^{\alpha_i}}{r_Np^{\alpha_N}}}\\
                                                  = &\frac{1}{r_Np^{\alpha_N}}\Big(1- \frac{\lambda+\sum_{i=1}^{N-1}r_ip^{\alpha_i}}{r_Np^{\alpha_N}}\Big) + \mathcal{O}(p^{-3\alpha_N+\alpha_{N-1}}) \\
                                                  = & \frac{r_Np^{\alpha_N}-\sum_{i=1}^{N-1}r_ip^{\alpha_i}}{(r_Np^{\alpha_N})^2} + \mathcal{O}(p^{-2\alpha_N}).
\end{align*}
The definition of the eigenvalue problem $A\varphi_n=\lambda_n\varphi_n$, integration by parts twice and the completeness of the eigenfunctions $(\varphi_n)_{n=1}^\infty$ in $L^2(\Omega)$ lead to
\begin{align}
  & \sum_{n=1}^\infty \lambda_n\langle u_0,\varphi_n\rangle\varphi_n (x) = \sum_{n=1}^\infty \langle u_0,\lambda_n\varphi_n\rangle\varphi_n (x) \nonumber\\
  = &  \sum_{n=1}^\infty \langle u_0,A\varphi_n\rangle\varphi_n (x) = \sum_{n=1}^\infty \langle Au_0,\varphi_n\rangle\varphi_n (x) = Au_0(x).\label{eqn:Au}
\end{align}
Substituting the last two identities, we obtain
\begin{align*}
  \mathcal{R}^*\hat u(x_0,p)
                & = p^{-1}\sum_{n=1}^\infty\langle u_0,\varphi_n\rangle\mathcal{R}^*\varphi_n(x_0) - \frac{r_Np^{\alpha_N}-\sum_{i=1}^{N-1}r_ip^{\alpha_i}}{p(r_Np^{\alpha_N})^2}  \mathcal{R}^*Au_0(x_0) + \mathcal{O}(p^{-2\alpha_N-1}).
\end{align*}
Now the Karamata-Feller Tauberian theorem \cite[Section XIII. 5]{Feller:1970} and the standard inverse Laplace transform relation
$\mathcal{L}^{-1} [p^{-s-1}] = \frac{t^{s}}{\Gamma(s)}$ for $s\geq0$
directly imply the first assertion. Similarly, the proof of Theorem \ref{thm:main1} and repeating the argument
lead to the following expansion as $p\to\infty$
\begin{align*}
  &\quad \mathcal{R}^*\hat u(x_0,p)  = \sum_{n=1}^\infty \frac{\hat \sigma(p) }{\lambda_n+\sum_{i=1}^{N}r_ip^{\alpha_i}}\langle f,\varphi_n\rangle\mathcal{R}^*\varphi_n(x_0)\\
                             & = \Big(\frac{1}{r_Np^{\alpha_N+a+1}}-\sum_{i=1}^{N-1}\frac{r_i}{r_N^2p^{2\alpha_N-\alpha_i+a+1}}\Big)
                             \sum_{n=1}^\infty\langle f,\varphi_n\rangle\mathcal{R}^*\varphi_n(x_0) + \mathcal{O}(p^{-2\alpha_N-a-1})\\
                             & = \Big(\frac{1}{r_Np^{\alpha_N+a+1}}-\sum_{i=1}^{N-1}\frac{r_i}{r_N^2p^{2\alpha_N-\alpha_i+a+1}}\Big)
                             \mathcal{R}^* f(x_0) + \mathcal{O}(p^{-2\alpha_N-a-1}).
\end{align*}
Then by the inverse Laplace transform and Karamata-Feller theorem \cite[Section XIII. 5]{Feller:1970}, we obtain the second expression.
\end{proof}

\begin{remark}\label{rmk:asymptotic}
The result in Proposition \ref{prop:asymptotic-u0} can also be seen as follows. Recall the
multinomial Mittag-Leffler function $E_{(\beta_1,\ldots,\beta_m),\beta_0}(z_1,\ldots,z_m)$
defined by \cite{LuchkoGorenflo:1999,Luchko:2011}
\begin{equation*}
  E_{(\beta_1,\ldots,\beta_m),\beta_0}(z_1,\ldots,z_m) = \sum_{k=0}^\infty \sum_{k_1+\ldots+k_m=k} \frac{(k;k_1,\ldots,k_m)\prod_{j=1}^mz_j^{k_j}}{\Gamma(\beta_0+\sum_{j=1}^m\beta_jk_j)},
\end{equation*}
where $0<\beta_0<2$, $0<\beta_j<1$ and $z_j\in\mathbb{C}$, $j=1,\ldots,m$, and the notation $(k;k_1,\ldots,k_m)$ denotes
the multinomial coefficient
\begin{equation*}
  (k;k_1,\ldots,k_m) = \frac{k!}{k_1!\ldots k_m!},\quad \mbox{for } k =k_1+\ldots+k_m,\quad k_1,\ldots,k_m\geq 0.
\end{equation*}
This function is a generalization of the classical two-parameter Mittag-Leffler function $E_{\alpha,\beta}(z)$, for $\alpha\in(0,2)$ and $\beta\in\mathbb{R}$ defined by $E_{\alpha,\beta}(z) = \sum_{k=0}^\infty \frac{z^k}{\Gamma(k\alpha+\beta)}$, for $z\in \mathbb{C}$,
which is an entire function generalizing the familiar exponential function \cite[Section 3.1]{Jin:2021}. Then assuming $r_N=1$, the solution $u$ to problem
\eqref{eqn:fde} with $ f\equiv0$ is given by (noting $0<\alpha_1<\ldots<\alpha_N<1$) \cite{Luchko:2011,LiLiuYamamoto:2015}
\begin{align*}
  u(x,t) &= \sum_{n=1}^\infty (1-\lambda_nt^{\alpha_N}E_{{\tilde{\boldsymbol\alpha}},1+\alpha_N}(-\lambda_nt^{\alpha_N},-r_1t^{\alpha_N-\alpha_1},\ldots,-r_{N-1}t^{\alpha_N-\alpha_{N-1}}))\langle u_0,\varphi_n\rangle\varphi_n(x),
\end{align*}
{with $\tilde{\boldsymbol{\alpha}}=(\alpha_N,\alpha_N-\alpha_1,\ldots,\alpha_N-\alpha_{N-1})$}
Thus, for small $t$, the identity \eqref{eqn:Au} implies
\begin{align*}
  \mathcal{R}^* u(x_0,t) = & \sum_{n=1}^\infty \langle u_0,\varphi_n\rangle \mathcal{R}^*\varphi_n(x_0) -\frac{t^{\alpha_N}}{\Gamma(\alpha_N+1)}\sum_{n=1}^\infty
   \lambda_n\langle u_0,\varphi_n\rangle\mathcal{R}^*\varphi_n(x_0) \nonumber\\
   &+\sum_{i=1}^{N-1}  \frac{r_it^{2\alpha_N-\alpha_i}}{\Gamma(2\alpha_N-\alpha_i+1)}\sum_{n=1}^\infty\lambda_n\langle u_0,\varphi_n\rangle\mathcal{R}^*\varphi_n(x_0) +\mathcal{O}(t^{2\alpha_N}) \nonumber\\
   =&\sum_{n=1}^\infty \mathcal{R}^*u_0(x_0)-\mathcal{R}^*Au_0(x_0)\Big(\frac{t^{\alpha_N}}{\Gamma(\alpha_N+1)}-\sum_{i=1}^{N-1}\frac{r_it^{2\alpha_N-\alpha_i}}{\Gamma(2\alpha_N-\alpha_i+1)}\Big) + \mathcal{O}(t^{2\alpha_N}).
\end{align*}
Thus we have deduced the desired asymptotic expansion in Proposition \ref{prop:asymptotic-u0} when $r_N=1$, and
the general case follows by a simple scaling argument. Note that the summation in the bracket
actually has the opposite sign of the leading term, and the constants are fully determined by $r_i$ and $\alpha_i$.
Meanwhile, by the Karamata-Feller Tauberian theorem  \cite[Section XIII. 5]{Feller:1970}, the condition on
the function $\sigma$ can be restated as $\sigma (t)= \frac{c_0}{\Gamma(a+1)}t^{a} + \mathcal{O}(t^{a+1})$ as $t\to 0 ^+$,
and the expression can be derived similarly using the identity
\begin{align*}
  &\frac{\d}{\d t}t^{\alpha_N} E_{\tilde{\boldsymbol{\alpha}},1+\alpha_N}
  (-\lambda_nt^{\alpha_N},-r_1t^{\alpha_N-\alpha_1},\cdots,-r_{N-1}t^{\alpha_N-\alpha_{N-1}})\\
  =& t^{\alpha_N-1}E_{\tilde{\boldsymbol{\alpha}},\alpha_N}
  (-\lambda_nt^{\alpha_N},-r_1t^{\alpha_N-\alpha_1},\cdots,-r_{N-1}t^{\alpha_N-\alpha_{N-1}}),
\end{align*}
and the fact that the solution $u(x,t)$ to problem  \eqref{eqn:fde} with $u_0\equiv0$ is given by
\begin{align*}
  u(t) &= \sum_{n=1}^\infty \int_0^ts^{\alpha_N-1}E_{\tilde{\boldsymbol{\alpha}},\alpha_N}(-\lambda_ns^{\alpha_N},
  -r_1s^{\alpha_N-\alpha_1},\ldots,-r_{N-1}s^{\alpha_N-\alpha_{N-1}}) \sigma(t-s)\d s \langle f,\varphi_n\rangle\varphi_n.
\end{align*}
\end{remark}

Motivated by Proposition \ref{prop:asymptotic-u0}, naturally one can develop
a procedure for numerically recovering the orders $\alpha_i$ and the corresponding weights $r_i$. The
most direct approach is to fit the fractional powers by a nonlinear least-squares formulation:
\begin{equation}\label{eqn:lsq}
  J(\mathbf{c},\boldsymbol{\beta}) = \tfrac12\|u(x_0,t)-f(\mathbf{c},\boldsymbol {\beta})\|_{L^2(0,T_0)}^2,
\end{equation}
with the regressor / model function $f_p(\mathbf{c},\boldsymbol \beta)$ given by
\begin{equation*}
f_p(\mathbf{c},\boldsymbol {\beta})=\left\{\begin{aligned}
c_0+\sum_{i=1}^{N}c_it^{\beta_i}, &\quad \mbox{if }u_0\not\equiv0, f\equiv0,\\
\sum_{i=1}^{N}c_it^{\beta_i}, & \quad \mbox{if } u_0\equiv0, f\not\equiv0,
\end{aligned}\right.
\end{equation*}
provided that the regressor $f(\mathbf{c},\boldsymbol\beta)$ approximates $u(x_0,t)$
well. Nonlinear least-squares problems of this type have been employed in statistics \cite{RoystonSauerbrei:2008}.
Once the parameters $\mathbf{c}$ and $\boldsymbol \beta$ are determined, the orders $\alpha_i$ and the
weights $r_i$ can be easily deduced. Numerical experiments indicate that this condition indeed holds under certain
assumptions: either $\alpha_N$ is close to one, or $T_0$ is sufficiently close to zero, which are however
not known a priori and potentially can restrict the range of application. This restriction is
severe when $\alpha_N$ is close to zero, and probably alternative models are needed. We shall discuss one alternative
based on rational approximation.

\subsection{Numerical results}
Now we present numerical results to illustrate the feasibility of recovering the orders and weights from time trace
data $g(t)=u(x_0,t)$, when the direct problem \eqref{eqn:fde} is equipped with a zero Neumann boundary
condition and the medium is unknown. Problem \eqref{eqn:lsq} uses 100 points uniformly
distributed on the interval $[0,T_0]$, for the time horizon $T_0$. It is minimized via the stand-alone
algorithm L-BFGS-B \cite{ByrdLuNocedal:1995}, implemented via
\texttt{MATLAB} wrapper from \url{https://www.mathworks.com/matlabcentral/fileexchange/35104-lbfgsb-l-bfgs-b-mex-wrapper}
(retrieved on May 6, 2021). L-BFGS-B is a quasi-Newton type algorithm that can take care of box constraints on the unknown,
and requires computing the value and gradient of the objective function. For the examples presented below, the
algorithm converges within tens of iterations (which of course depends strongly on the initial guess), and thus
the overall procedure is fairly efficient. Note that problem
\eqref{eqn:lsq} is highly nonconvex for both model functions (i.e., fractional polynomials and rational functions),
and thus a good initial guess is needed in order to ensure reasonable recovery.

First we discuss a simple yet illuminating example.
\begin{example}\label{exam:single}
The domain $\Omega=(0,1)$, $Au=-u''+u$ with a zero Neumann boundary condition, $u_0=\cos\pi x$, and $f\equiv0$.
The measurement point $x_0$ is the left end point $x_0=0$.
{\rm (i)} One single-order with $\alpha\in (0,1)$ and {\rm (ii)} Two terms with $0<\alpha_1<\alpha_2<1$ and the corresponding weights $r_1$ and $r_2=1$.
Note that $u_0$ is actually an eigenfunction of the operator $A$, with the corresponding eigenvalue $\lambda = \pi^2+1$,
and thus the direct problem essentially reduces to an ODE. The goal is to recover $\boldsymbol
\alpha$ and $\lambda$ from the data $g(t)=u(x_0,t)$.
\end{example}

Note that for this example, the data $g(t)$ admits a closed form
\begin{equation*}
  g(t) =  \left\{\begin{aligned}
      E_{\alpha,1}(-\lambda t^\alpha), & \quad \mbox{case (i)},\\
      1-\lambda t^{\alpha_2}E_{(\alpha_2,\alpha_2-\alpha_1),1+\alpha_2}(-\lambda t^{\alpha_2},-r_1t^{\alpha_2-\alpha_1}), &\quad \mbox{case (ii)}.
  \end{aligned}\right.
\end{equation*}
By Proposition \ref{prop:asymptotic-u0}, the asymptotic of $g(t)$ as $t\to 0^+$ is given by
\begin{equation*}
   g(t) = \left\{\begin{aligned}
      1-\lambda \frac{t^\alpha}{\Gamma(\alpha+1)} + \mathcal{O}(t^{2\alpha}), & \quad \mbox{case (i)},\\
      1-\lambda \frac{t^{\alpha_2}}{\Gamma(\alpha_2+1)} + \lambda \frac{r_1t^{2\alpha_2-\alpha_1}}{\Gamma(2\alpha_2-\alpha_1+1)} +\mathcal{O}(t^{2\alpha_2}), &\quad \mbox{case (ii)}.
  \end{aligned}\right.
\end{equation*}

\begin{figure}[hbt!]
\centering
\setlength{\tabcolsep}{0pt}
\begin{tabular}{cccc}
\includegraphics[height=0.20\textwidth]{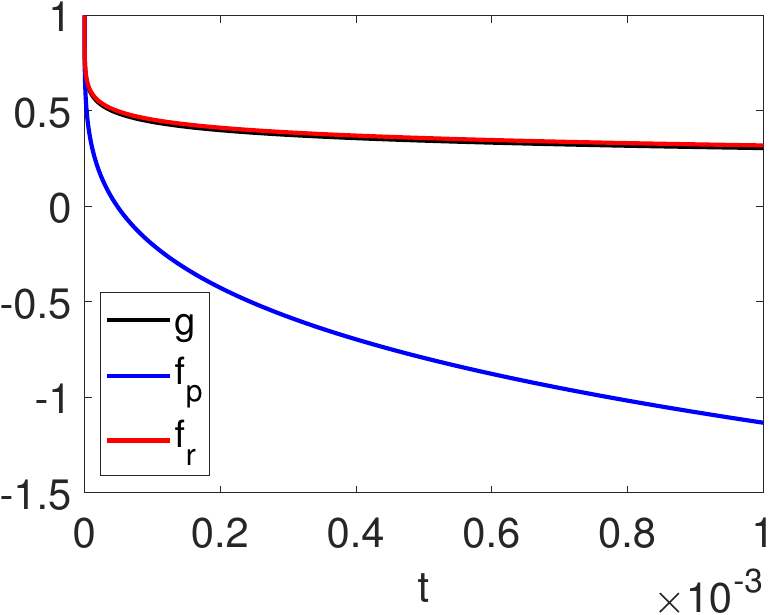} & \includegraphics[height=0.20\textwidth]{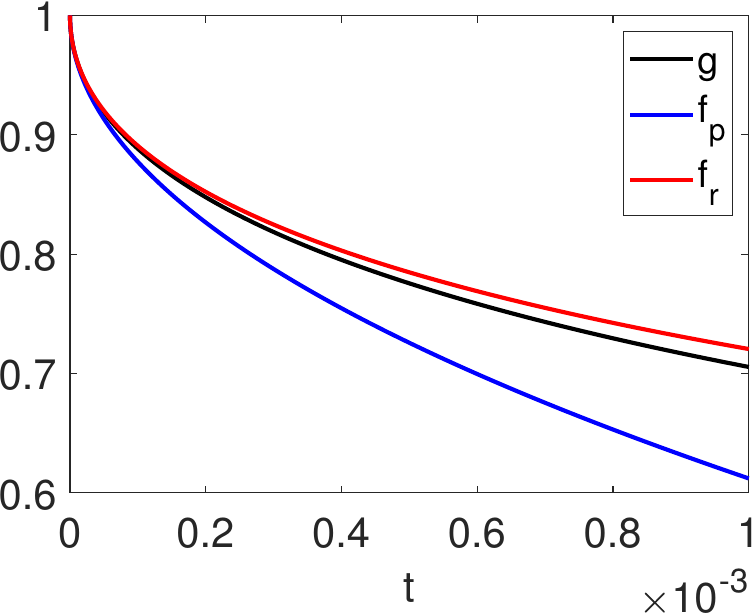}
&\includegraphics[height=0.20\textwidth]{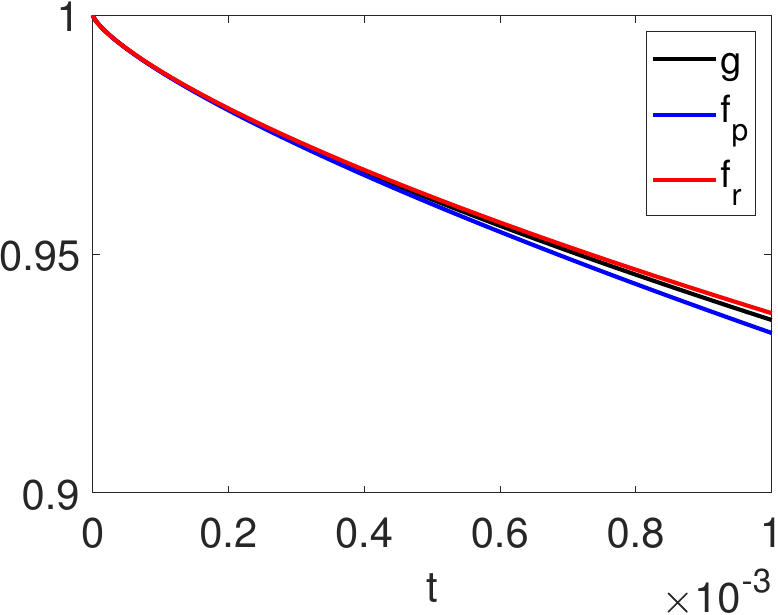} & \includegraphics[height=0.20\textwidth]{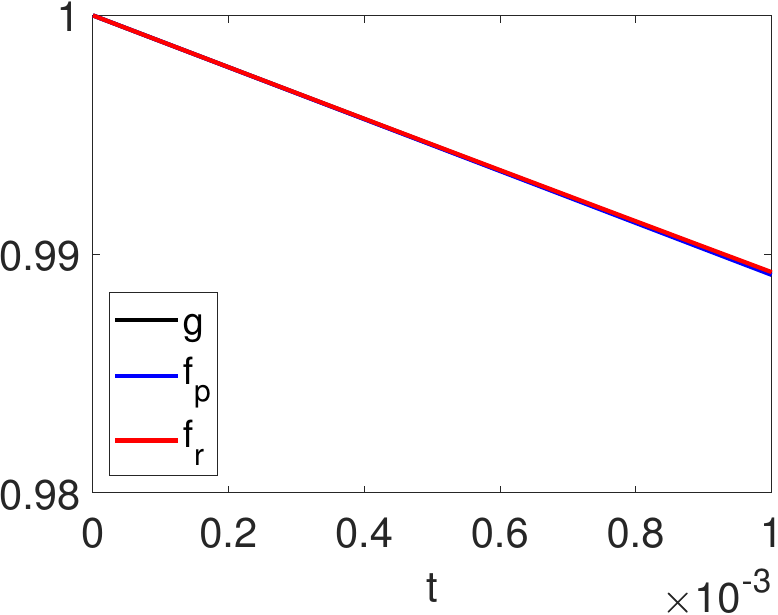}\\
(a) $\alpha=0.25$ & (b) $\alpha=0.5$ & (c) $\alpha=0.75$ & (d) $\alpha=1.00$
\end{tabular}
\caption{The Mittag-Leffler function $g(t)=E_{\alpha,1}(-\lambda t^\alpha)$ and the fractional polynomial approximation $f_p$ and rational
approximation $f_r$ for Example \ref{exam:single}(i).}\label{fig:mlf}
\end{figure}

\begin{figure}[hbt!]
\centering
\setlength{\tabcolsep}{0pt}
\begin{tabular}{cccc}
\includegraphics[height=0.20\textwidth]{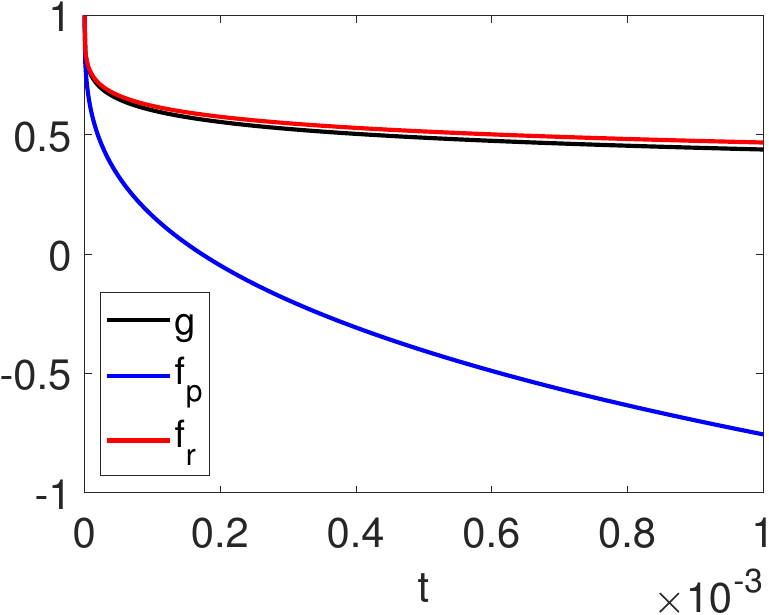} & \includegraphics[height=0.20\textwidth]{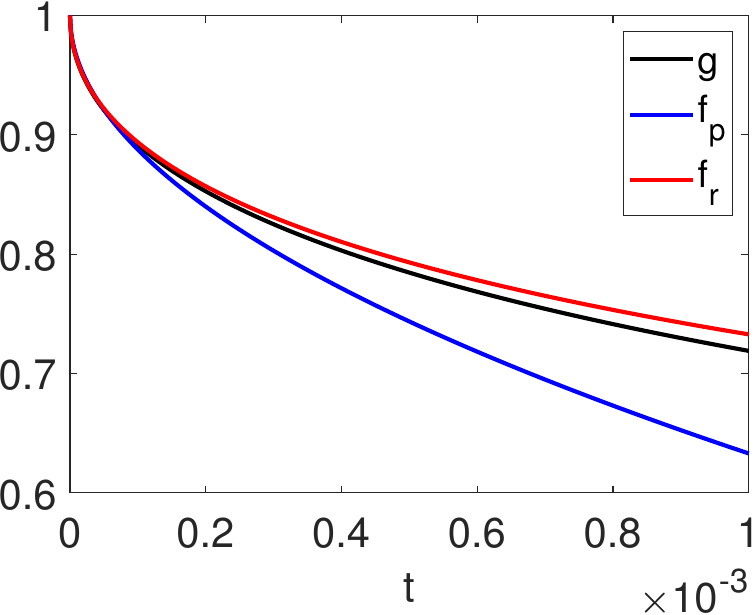}
&\includegraphics[height=0.20\textwidth]{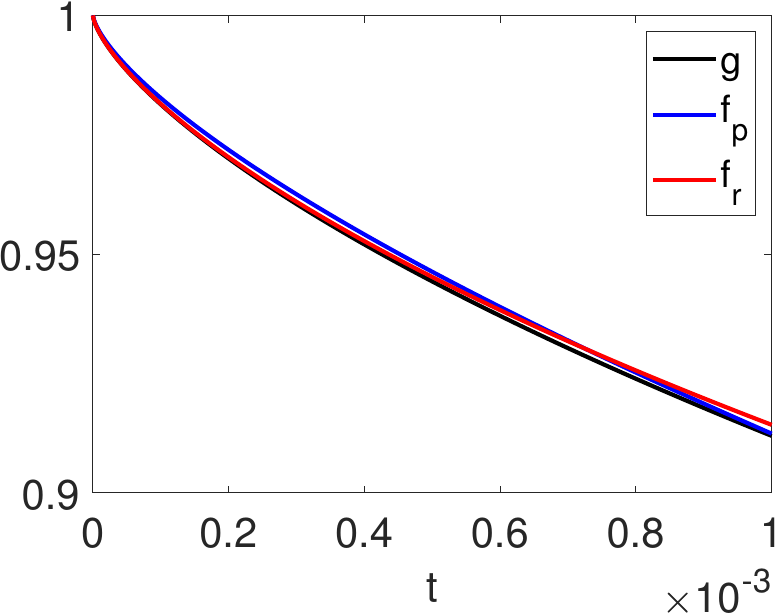} & \includegraphics[height=0.20\textwidth]{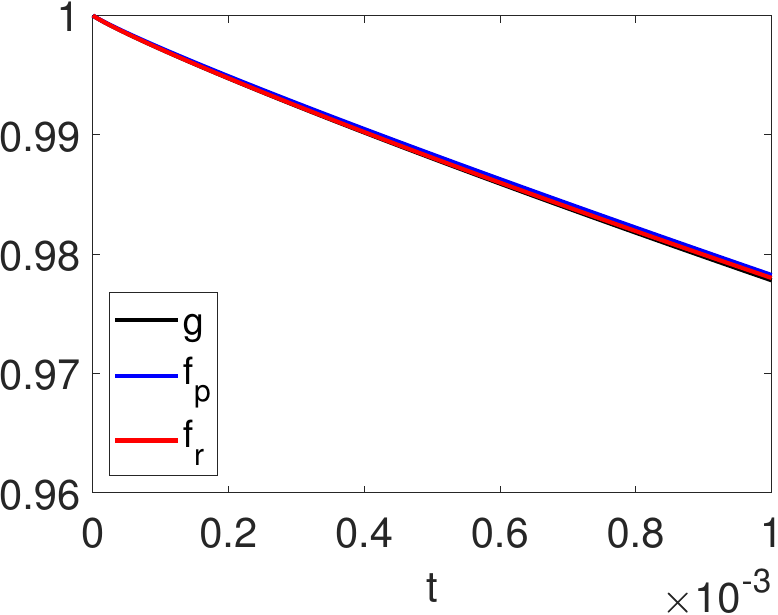}\\
(a) $\boldsymbol{\alpha}=(0.2,0.3)$ & (b) $\boldsymbol{\alpha}=(0.2,0.5)$ & (c) $\boldsymbol{\alpha}=(0.2,0.7)$ & (d) $\boldsymbol{\alpha}=(0.2,0.9)$
\end{tabular}
\caption{The multinomial Mittag-Leffler function $g(t)=1-\lambda t^{\alpha_2}E_{(\alpha_2,\alpha_2-\alpha_1),
1+\alpha_2}(-\lambda t^{\alpha_2},-r_1t^{\alpha_2-\alpha_1})$ and the fractional polynomial approximation $f_p$ and rational
approximation $f_r$ for Example \ref{exam:single}(ii) with $r_1=0.5$.}\label{fig:mml}
\end{figure}

Since the feasibility of the formulation \eqref{eqn:lsq} resides on the accuracy of the fractional
polynomial approximation, denoted by $f_p(\mathbf{c},\boldsymbol\alpha)$ below, we investigate the
accuracy of the approximation. We evaluate the function $E_{\alpha,1}(-\lambda t^\alpha)$ over
a small time interval $[0,T_0]$ by an algorithm from \cite{SeyboldHilfer:2008}, and the multinomial
Mittag-Leffler function by summing the power series over $k$ (truncated at $k=100$), since the series
converges rapidly for a small argument (we are not aware of any algorithm for evaluating the multinomial
Mittag-Leffler functions). The numerical results are shown in Figs. \ref{fig:mlf}
and \ref{fig:mml}, for the single- and two-term, respectively. Note that the function $f_p$ can
approximate the target function $g(t)$ over a suitable interval $[0,T_0]$ accurately, but the size $T_0$ depends strongly
on the order $\alpha$. In particular, as the order $\alpha$ tends to zero, the size $T_0$ should
shrink to zero so as to maintain the desired accuracy. Thus, the least-squares approach \eqref{eqn:lsq}
based on the model $f_p$ should only make use of data within $[0,T_0]$ with a tiny $T_0$, in order to
handle a broad range of $\alpha$. This observation remains largely valid in the two-term cases, i.e.,
the $T_0$ should be chosen small when both orders are small, cf. Fig. \ref{fig:mml}.

These empirical observations naturally motivate the question whether there is actually a better model
than fractional polynomials for the numerical recovery. In the simple setting, the answer is affirmative. One promising model is
the lowest-order rational approximation $f_r$ given by
\begin{equation*}
   f_r(t) = \left\{\begin{aligned}
     \frac{1}{1+\lambda\frac{t^\alpha}{\Gamma(\alpha+1)}}, & \quad\mbox{case (i)},\\
     \frac{1}{1+\lambda\frac{t^{\alpha_2}}{\Gamma(\alpha_2+1)}-\lambda\frac{r_1t^{2\alpha_2-\alpha_1}}{\Gamma(2\alpha_2-\alpha_1+1)}}, &\quad \mbox{case (ii)}.
   \end{aligned}\right.
\end{equation*}
Actually, it is known that for (i), $f_r$ is actually an upper bound on $E_{\alpha,1}(-\lambda t^\alpha)$ (see e.g.,
\cite{Simon:2014} or \cite[Theorem 3.6]{Jin:2021}, and empirically observed by Mainardi \cite{Mainardi:2014}).
Numerically, the model $f_r$ is asymptotically tight as $t\to 0^+$, just as the model $f_p$, but it
approximates better $g(t)$. We are not aware of the two-term case in the existing literature, but one can glim a hint
on the approximation from \cite[Lemma 6.1]{VergaraZacher:2015} and fractional polynomial approximation of the
resolvent kernel. The numerical results in Figs. \ref{fig:mlf} and \ref{fig:mml} show clearly that the model $f_r$
is indeed a more accurate approximation to the data $g(t)$, for both single- and two-term cases, although there
is still no rigorous proof of the empirical observation yet. Thus, it is conceived that the rational model $f_r$ is
better suited as the regressor for the numerical recovery via the nonlinear least-squares approach \eqref{eqn:lsq}
(of course, modulus the challenge of nonlinear optimization).

Next we present results for numerical recovery of the orders for Example \ref{exam:single}. The numerical results
are presented in Tables \ref{tab:single-i} and \ref{tab:single-ii} for cases (i) and (ii), respectively.
It is observed that for a fixed $\alpha=0.7$, the recovered $\alpha$ represents a very accurate
approximation to the exact one, when the time horizon $T_0$ is sufficiently close to zero, and then
both models $f_p$ and $f_r$ have comparable accuracy. However, when $T_0$ increases, the results by
the model $f_r$ is much more accurate than that by the model $f_p$, which agrees with the preceding
observation that $f_r$ is more broadly valid as an approximation to the target function $g(t)$. Meanwhile,
when $T_0$ is fixed at 1e-6, the accuracy of the recovered order $\alpha$ deteriorates steadily as the
true $\alpha$ decreases towards zero, for either model function, although the results by the model function $f_r$ are far
more accurate, especially when the exact $\alpha$ is small. This can be attributed to the
behavior of $E_{\alpha,1}(-\lambda t^\alpha)$: at small $\alpha$, $E_{\alpha,1}(-\lambda t^\alpha)$ reaches a
quasi-steady state very rapidly, cf. Fig. \ref{fig:mlf}(a), and the fractional polynomial model $f_p$
fails to accurately capture the behavior, whereas the rational model $f_r$ does so more closely.

\begin{table}[hbt!]
\centering
\caption{Numerical results for Example \ref{exam:single}(i), $\lambda=\pi^2+1\approx1.0870$, and the algorithm L-BFGS-B
is always initialized to $\alpha=0.5$, and $c$ by the standard linear least-squares method. (a) Fixed $\alpha=0.7$,
  recovered with different $T_0$, and (b) recovered with $T_0=$1e-6, for different $\alpha$.\label{tab:single-i}}
\begin{threeparttable}
  \subfloat[]{
  \begin{tabular}{ccccc}
  \toprule
  \multicolumn{1}{c}{} & \multicolumn{2}{c}{$f_p$} &\multicolumn{2}{c}{$f_r$}\\
  \cmidrule(lr){2-3} \cmidrule(lr){4-5}
  $T_0$ & $\alpha$ & $\lambda$ & $\alpha$ & $\lambda$\\
  \midrule
  1e-7 &  0.6993 &  1.074e1 &  0.6995 &  1.078e1\\
  1e-6 &  0.6998 &  1.083e1 &  0.7001 &  1.088e1\\
  1e-5 &  0.6989 &  1.071e1 &  0.7005 &  1.095e1\\
  1e-4 &  0.6947 &  1.022e1 &  0.7026 &  1.121e1\\
  1e-3 &  0.6742 &  8.532e0 &  0.7129 &  1.228e1\\
  1e-2 &  0.5856 &  4.779e0 &  0.7570 &  1.648e1\\
  1e-1 &  0.3455 &  1.660e0 &  0.8563 &  2.690e1\\
  \bottomrule
  \end{tabular}}
  \subfloat[]{
  \begin{tabular}{ccccc}
  \toprule
  \multicolumn{1}{c}{} & \multicolumn{2}{c}{$f_p$} &\multicolumn{2}{c}{$f_r$}\\
  \cmidrule(lr){2-3} \cmidrule(lr){4-5}
  $\alpha$ & $\alpha$ & $\lambda$ & $\alpha$ & $\lambda$\\
  \midrule
     0.30 &  0.2669 & 5.916e0 & 0.3034 & 1.157e1\\
     0.40 &  0.3891 & 9.008e0 & 0.4019 & 1.124e1\\
     0.50 &  0.4969 & 1.032e1 & 0.5008 & 1.102e1\\
     0.60 &  0.5991 & 1.072e1 & 0.6003 & 1.092e1\\
     0.70 &  0.6998 & 1.083e1 & 0.7001 & 1.088e1\\
     0.80 &  0.7995 & 1.079e1 & 0.7990 & 1.072e1\\
     0.90 &  0.8937 & 9.906e0 & 0.8909 & 9.521e0\\
  \bottomrule
  \end{tabular}}
  \end{threeparttable}
\end{table}

\begin{table}[hbt!]
\centering
\caption{Numerical results for Example \ref{exam:single}(ii), $\lambda
=\pi^2+1\approx1.0870$, and the algorithm L-BFGS-B, and the corresponding
$c$ by the standard linear least-squares method.\label{tab:single-ii}}
\begin{threeparttable}
  \subfloat[$\boldsymbol{\alpha}=(0.6,0.9)$, L-BFGS-B initialized to $\boldsymbol\alpha=(0.4,0.8)$]{
  \begin{tabular}{ccccccccc}
  \toprule
  \multicolumn{1}{c}{} & \multicolumn{4}{c}{$f_p$} &\multicolumn{4}{c}{$f_r$}\\
  \cmidrule(lr){2-5} \cmidrule(lr){6-9}
  $T_0$ & $\alpha_1$ & $\alpha_2$ & $\lambda$ & $r_1$ & $\alpha_1$ & $\alpha_2$ & $\lambda$ & $r_1$ \\
  \midrule
  1e-7 & 0.5971 &  0.8985 &  1.056e1 &        0  &   0.5971 &  0.8985 &  1.056e1 &        0  \\
  1e-6 & 0.5943 &  0.8972 &  1.032e1 &  5.49e-12 &   0.5944 &  0.8972 &  1.033e1 &  6.61e-13 \\
  1e-5 & 0.5887 &  0.8943 &  9.938e0 &  1.04e-11 &   0.5890 &  0.8945 &  9.963e0 &  1.53e-11 \\
  1e-4 & 0.5766 &  0.8883 &  9.284e0 &  1.59e-9  &   0.5794 &  0.8897 &  9.432e0 &  1.25e-09 \\
  1e-3 & 0.4227 &  0.8955 &  1.015e1 &  2.08e0   &   0.5844 &  0.8922 &  9.830e0 &  3.36e-01 \\
  1e-2 & 0.3530 &  0.8931 &  9.860e0 &  2.82e0   &   0.5788 &  0.8894 &  9.220e0 &        0  \\ \bottomrule
  \end{tabular}}

  \subfloat[$\boldsymbol\alpha=(0.5,0.7)$, L-BFGS-B initialized to $\boldsymbol\alpha=(0.2,0.6)$]{
  \begin{tabular}{ccccccccc}
  \toprule
  \multicolumn{1}{c}{} & \multicolumn{4}{c}{$f_p$} &\multicolumn{4}{c}{$f_r$}\\
  \cmidrule(lr){2-5} \cmidrule(lr){6-9}
  $T_0$ & $\alpha_1$ & $\alpha_2$ & $\lambda$ & $r_1$ & $\alpha_1$ & $\alpha_2$ & $\lambda$ & $r_1$ \\
  \midrule
  1e-7 & 0.3889 &  0.6944 &  9.616e0 &  2.28e-12 &   0.3890 &  0.6945 &  9.626e0 &  5.67e-11\\
  1e-6 & 0.3824 &  0.6912 &  9.129e0 &  2.41e-11 &   0.3830 &  0.6915 &  9.174e0 &  5.53e-11\\
  1e-5 & 0.3713 &  0.6856 &  8.466e0 &  1.03e-9 &   0.3742 &  0.6871 &  8.639e0 &  8.70e-10\\
  1e-4 & 0.3914 &  0.6959 &  9.954e0 &  1.50e0   &   0.5870 &  0.7010 &  1.160e1 &  4.71e-1 \\
  1e-3 & 0.2804 &  0.6906 &  9.156e0 &  2.60e0   &   0.3625 &  0.6803 &  8.051e0 &  9.89e-2 \\
  1e-2 & 0.5378 &  0.7315 &  1.633e1 &  1.41e0   &  -0.5278 &  0.6860 &  8.261e0 &        0 \\ \bottomrule
  \end{tabular}}
  \end{threeparttable}
\end{table}

The numerical results for the two term case in Example \ref{exam:single}(ii) are summarized in
Table \ref{tab:single-ii}. It is observed that with proper initialization, the method does
recover the orders $(\alpha_1,\alpha_2)$ to a reasonable accuracy, for a wide range of
the time horizon $T_0$, with the accuracy of $\alpha_2$ being higher than that of $\alpha_1$. The
latter can be attributed to the asymptotic expansion: the term involving $\alpha_2$ is dominating
in the expansion and much easier to estimate than the remaining terms. Indeed, if $T_0$ is
sufficiently small, all other terms are essentially negligible, comparing the results in
Figs. \ref{fig:mlf} and \ref{fig:mml}, and consequently, the order $\alpha_1$ and the weight
$r_1$ cannot be estimated reliably at all. This can also been seen from the following slightly
more refined expansion in fractional polynomials:
\begin{align*}
   1-\lambda t^{\alpha_2}E_{(\alpha_2,\alpha_2-\alpha_1),1+\alpha_2}(-\lambda t^{\alpha_2},-r_1t^{\alpha_2-\alpha_1})
 = f_p(t) + \lambda^2\frac{t^{2\alpha_2}}{\Gamma(2\alpha_2+1)}+\mathcal{O}(t^{3\alpha_2-\alpha_1}).
\end{align*}
This clearly shows the potential pitfalls in the order recovery: the next term can
have comparable magnitude with the last term $f_p$ when $\alpha_1$ is close to zero, and
the nonlinear procedure attempts to approximate it with the leading terms, thereby
significantly affecting the recovery accuracy.
The accuracy of the recovered $\lambda$ is also reasonable, except for fairly large $T_0$.
However, the accuracy of $r_1$ is poor in all cases, due to the aforementioned reasons. Also
as the orders decrease, it is becoming increasingly more challenging for the numerical recovery, which
agrees with the empirical observation from Example \ref{exam:single}(i), cf. Table \ref{tab:single-ii}(b).
This is attributed to the rapid decay near $t=0$ so that the model functions are not accurate.
These results partly confirm the assertion in Theorem \ref{thm:main2}: the recovery is indeed
possible, however, numerically this can still be a big challenge, depending on the magnitude
of the sought-for orders. It is of great interest to develop further remedies to tackle the numerical issues.

The next example illustrates the feasibility of the approach in the general setting.
\begin{example}\label{exam:mix}
The domain $\Omega=(0,1)^2$, $Au=-\Delta u$ with zero Neumann boundary condition, with a two-term model
with $0<\alpha_1<\alpha_2<1$ and the corresponding weights $r_1$ and $r_2=1$. The measurement point $x_0$
is the {vertex point} $x_0=(0,0)$. {Consider the following two cases:}
{\rm (i)} $u_0(x_1,x_2)=\cos\pi x_1 \cos \pi x_2 + \frac{1}{4}(\cos 2\pi x_1\cos \pi x_2+\cos\pi x_1\cos2\pi x_2)+\frac{1}{8}\cos 2\pi x_1\cos 2\pi x_2$ and $f\equiv0$ and {\rm (ii)} $u_0\equiv0$ and $f(x_1,x_2,t)= \cos\pi x_1 \cos \pi x_2 + \frac{1}{2}(\cos 2\pi x_1\cos \pi x_2+\cos\pi x_1\cos2\pi x_2)+\frac{1}{4}(\cos 3\pi x_1\cos \pi x_2+\cos\pi x_1\cos 3\pi x_2)$.
The goal is to recover the following quantities: $\boldsymbol \alpha$, $r_1$ and $Au(x_0)$ or $f(x_0)$ for {cases {\rm(i)} and {\rm(ii)}}, respectively, from the data $g(t)=u(x_0,t)$.
\end{example}

Note that for this example, the data $g(t)=u(x_0,t)$ can be generated using the standard Galerkin finite element method in
space and convolution quadrature in time \cite{JinLazarovZhou:2016sisc}. Below we employ series expansion
using the multinomial Mittag-Leffler function, cf. Remark \ref{rmk:asymptotic}, so as to minimize
the discretization errors, and like before, it is evaluated by means of series summation
(truncated at $k=100$). In case (ii), the model functions $f_p$ and $f_r$ are given respectively by
\begin{align*}
  f_p(t) &= f(x_0)\Big(\frac{t^{\alpha_2}}{\Gamma(\alpha_2+1)}-\frac{r_1t^{2\alpha_2-\alpha_1}}{\Gamma(2\alpha_2-\alpha_1+1)}\Big),\\
  f_r(t) &= f(x_0)\Big(1-\frac{1}{1+\frac{t^{\alpha_2}}{\Gamma(\alpha_2+1)}-\frac{r_1t^{2\alpha_2-\alpha_1}}{\Gamma(2\alpha_2-\alpha_1+1)}}\Big),
\end{align*}
where the latter follows by direct analogy. The numerical results are given in Table \ref{tab:mix}. The observations
from Example \ref{exam:single} remain largely valid. In both cases (i) and (ii), the orders can be accurately recovered,
provided that the time horizon $T_0$ is sufficiently small (so that the model functions are accurate
approximations). The accuracy of the two models are largely comparable to each other in either case, and deteriorates
steadily as $T_0$ increases. Also the accuracy of the recovered $Au(x_0)=$5.428e1 and $f(x_0)=2.500$ is fair. Just as
expected, the accuracy of the estimated $r_1$ is poor in all cases, irrespective of the time horizon $T_0$, which also
agrees with preceding observations. These numerical experiments not only confirm the possibility of uniquely recovery as indicated
by Theorems \ref{thm:main1} and \ref{thm:main2}, but also illustrate the pitfalls in developing practical recovery schemes.

\begin{table}[hbt!]
\centering
\caption{Numerical results for Example \ref{exam:mix}, and the algorithm L-BFGS-B, and the corresponding
$c$ by the standard linear least-squares method.\label{tab:mix}}
\begin{threeparttable}
  \subfloat[$\boldsymbol{\alpha}=(0.5,0.8)$, L-BFGS-B initialized to $\boldsymbol\alpha=(0.3,0.7)$]{
  \begin{tabular}{ccccccccc}
  \toprule
  \multicolumn{1}{c}{} & \multicolumn{4}{c}{$f_p$} &\multicolumn{4}{c}{$f_r$}\\
  \cmidrule(lr){2-5} \cmidrule(lr){6-9}
  $T_0$ & $\alpha_1$ & $\alpha_2$ & $Au(x_0)$ & $r_1$ & $\alpha_1$ & $\alpha_2$ & $Au(x_0)$ & $r_1$ \\
  \midrule
  1e-8 & 0.4992 & 0.7996 & 5.380e1 &  0 & 0.4992 & 0.7996 & 5.381e1 & 0 \\
  1e-7 & 0.4984 & 0.7992 & 5.340e1 &  0 & 0.4985 & 0.7992 & 5.344e1 & 0 \\
  1e-6 & 0.4965 & 0.7983 & 5.261e1 &  0 & 0.4971 & 0.7985 & 5.284e1 & 0 \\
  1e-5 & 0.4913 & 0.7956 & 5.079e1 &  0 & 0.4946 & 0.7973 & 5.196e1 & 0 \\
  1e-4 & 0.4716 & 0.7858 & 4.551e1 &  0 & 0.4920 & 0.7960 & 5.118e1 & 0 \\
  \bottomrule
  \end{tabular}}

  \subfloat[$\boldsymbol{\alpha}=(0.5,0.7)$, L-BFGS-B initialized to $\boldsymbol\alpha=(0.3,0.6)$]{
  \begin{tabular}{ccccccccc}
  \toprule
  \multicolumn{1}{c}{} & \multicolumn{4}{c}{$f_p$} &\multicolumn{4}{c}{$f_r$}\\
  \cmidrule(lr){2-5} \cmidrule(lr){6-9}
  $T_0$ & $\alpha_1$ & $\alpha_2$ & $f(x_0)$ & $r_1$ & $\alpha_1$ & $\alpha_2$ & $f(x_0)$ & $r_1$ \\
  \midrule
  1e-8 &    0.4964 & 0.6982 & 2.393 &        0  & 0.4964 & 0.6982 & 2.393 & 1.29e-10\\
  1e-7 &    0.4941 & 0.6970 & 2.343 & 1.27e-10  & 0.4941 & 0.6970 & 2.343 & 9.13e-10\\
  1e-6 &    0.5023 & 0.7015 & 2.583 & 6.41e-1   & 0.4896 & 0.6948 & 2.261 & 2.11e-9 \\
  1e-5 &    0.2083 & 0.6974 & 2.366 & 8.56e0    & 0.2167 & 0.6975 & 2.367 & 7.54e0  \\
  1e-4 &    0.1669 & 0.6969 & 2.343 & 1.23e1    & 0.1668 & 0.6967 & 2.337 & 1.16e1  \\
  \bottomrule
  \end{tabular}}
  \end{threeparttable}
\end{table}

\section{Conclusions}

In this work we have proved the unique recovery of multiple fractional orders and the associated
weights in a multi-term time-fractional diffusion model from the observation at one point on the boundary,
based on the asymptotics of the solution at small time and the time analyticity of the solution. We have
also discussed the numerical recovery based on asymptotic expansion / rational approximation, and
demonstrated the feasibility of the least-squares approach for recovering the highest
order and the weight.

It is of much interest to study related inverse problems for more complex  anomalous
diffusion models, e.g., distributed-order or variable-order. It is unclear whether
the one-point observation is sufficient for the unique recovery of the order distribution in these models, but a
partial determination, e.g., support in the distributed order, might be possible. Further,
it remains an outstanding challenge to develop stable and accurate numerical procedures for
recovering all fractional orders, by properly overcoming the difficulty with unknown media.

\bibliographystyle{abbrv}
\bibliography{frac}

\end{document}